\documentclass[10pt,a4paper]{amsart}
\usepackage[all]{xy}

\newtheorem{thm}{Theorem}[section]
\newtheorem{prop}[thm]{Proposition}
\newtheorem{lem}[thm]{Lemma}

\newcounter{thmlc}
\newenvironment{thmlist}
{
    \setcounter{thmlc}{0}
    \begin{list}
        {\roman{thmlc})}
        {
            \usecounter{thmlc}
            \setlength{\labelsep}{1mm}
            \setlength{\labelwidth}{5mm}
            \setlength{\leftmargin}{6mm}
            \setlength{\topsep}{0mm}
        }
}
{
    \end{list}
}

    \DeclareMathOperator*{\bigtimes}{\times}
    
    \DeclareMathOperator{\cen}{C}

    \DeclareMathOperator{\der}{Der}

    \DeclareMathOperator{\ext}{Ext}
    \DeclareMathOperator{\gr}{gr}

    \DeclareMathOperator{\hc}{HC}
    \DeclareMathOperator{\hh}{HH}
    \DeclareMathOperator{\hmm}{Hom}
    \DeclareMathOperator{\id}{id}
    
    \DeclareMathOperator{\im}{Im}

    \DeclareMathOperator{\mor}{Mor}
    \DeclareMathOperator*{\osum}{\oplus}
    \DeclareMathOperator{\pr}{pr}

\newcommand{\acat}[1]{\mathbf{a}_{#1}}
\newcommand{\Acat}[1]{\mathbf{A}_{#1}}

\newcommand{\cacat}[1]{\mathbf{\hat{a}}_{#1}}
\newcommand{\cpd}{*}

\newcommand{\defm}{\textup{Def}_{\mathbf{M}}}
\newcommand{\dvd}{\mid}

\newcommand{\dx}{\partial x}
\newcommand{\dy}{\partial y}
\newcommand{\ff}{\textup{F}}
\newcommand{\fg}{\textup{G}}
\newcommand{\fff}{\hat{\textup{F}}}
\newcommand{\ffg}{\hat{\textup{G}}}

\newcommand{\h}[1]{\textup{h}_{#1}}
\newcommand{\infdef}{\textup{T}^1}

\newcommand{\mfam}{\mathbf{M}}

\newcommand{\ndv}{\not \mspace{5.0mu} \dvd}

\newcommand{\obstr}{\textup{o}(u,M_S)}
\newcommand{\obstrdef}{\textup{T}^2}
\newcommand{\ol}{\overline}
\newcommand{\ord}[1]{|#1|}

\newcommand{\sets}{\mathbf{Sets}}
\newcommand{\test}[2]{R(#1,#2)}
\newcommand{\tg}[1]{t_{#1}}
\newcommand{\ul}{\underline}

\begin{document}

\title[Noncommutative deformations of modules]{An introduction to noncommutative
deformations of modules}
\author{Eivind Eriksen}
\address{Institute of Mathematics, University of Warwick \\ Coventry CV4 7AL, UK}
\email{eriksen@maths.warwick.ac.uk}

\thanks{This research has been supported by a Marie Curie Fellowship
of the European Community programme "Improving Human Research
Potential and the Socio-economic Knowledge Base" under contract
number HPMF-CT-2000-01099}

\begin{abstract}
Let $k$ be an algebraically closed (commutative) field, let $A$ be
an associative $k$-algebra, and let $\mfam = \{ M_1, \dots, M_p
\}$ be a finite family of left $A$-modules. We study the
simultaneous formal deformations of this family, described by the
noncommutative deformation functor $\defm: \acat p \to \sets$
introduced in Laudal \cite{lau02}. In particular, we prove that
this deformation functor has a pro-representing hull, and describe
how to calculate this hull using the cohomology groups
$\ext^n_A(M_i,M_j)$ and their matric Massey products.
\end{abstract}

\maketitle

\section*{Introduction}

In this paper, I shall give an elementary introduction to the
noncommutative deformation theory for modules, due to Laudal. This
theory, which generalizes the classical deformation theory for
modules, was introduced by Laudal in \cite{lau02}. Earlier
versions of this material appeared in the preprints Laudal
\cite{lau95}, \cite{lau95-2}, \cite{lau96}, \cite{lau98-2},
\cite{lau00}.

This noncommutative deformation theory has several applications:
In the paper Laudal \cite{lau02}, Laudal used it to construct
algebras with a prescribed set of simple modules, and also to
study the moduli space of iterated extensions of modules. In the
preprint Laudal \cite{lau00}, he also showed that this theory is a
useful tool in the study of algebras, and in establishing a
noncommutative algebraic geometry.

These applications are an important part of the motivation for the
noncommutative deformation theory. But we shall not go into the
details of these applications in this elementary introduction.
Instead, we refer to the papers and preprints of Laudal mentioned
above for applications and further developments of the theory.

Throughout this paper, we shall fix the following notations: Let
$k$ be an algebraically closed (commutative) field, let $A$ be an
associative $k$-algebra, and let $\mfam = \{ M_1, \dots, M_p \}$
be a finite family of left $A$-modules. Notice that this notation
differs from Laudal's: While Laudal considers families of right
modules in all his paper, I consider families of left modules. Of
course, the difference is only in the appearance --- the resulting
theories are obviously equivalent.

We shall present a noncommutative deformation functor $\defm:
\acat p \to \sets$, which describes the simultaneous formal
deformations of the family $\mfam$ of left $A$-modules.
Furthermore, we shall prove that this deformation functor has a
pro-representable hull $(H,\xi)$ when the family $\mfam$ satisfy a
certain finiteness condition. We shall also describe a method for
finding the pro-representable hull explicitly.

In section \ref{s:cat}, we describe the category $\acat p$. It is
a full sub-category of the category $\Acat p$ of $p$-pointed
$k$-algebras. The objects of $\Acat p$ are the $k$-algebras $R$
equipped with $k$-algebra homomorphisms $k^p \to R \to k^p$, such
that the composition $k^p \to k^p$ is the identify. For any such
object, $R = (R_{ij})$ is a $k$-algebra of $p \times p$ matrices.
The radical of this object is the ideal $I(R) = \ker(R \to k^p)
\subseteq R$. The category $\acat p$ is the full sub-category of
$\Acat p$ consisting of objects such that $R$ is Artinian and
complete in the $I(R)$-adic topology.

In section \ref{s:ncdefm}, we describe the noncommutative
deformation functor associated to the family $\mfam$ of left
$A$-modules,
    \[ \defm: \acat p \to \sets \]
It is constructed in the following way: Let $R$ be an object of
$\acat p$, and consider the vector space $M_R = (M_i \otimes_k
R_{ij})$, equipped with the natural right $R$-module structure
induced by the multiplication in $R$. A deformation of $\mfam$ to
$R$ consists of the following data:
\begin{itemize}
\item A left $A$-module structure on $M_R$ making $M_R$ a left $A
\otimes_k R^{\text{op}}$-module,
\item Isomorphisms $\eta_i: M_R \otimes_R k_i \to M_i$ of left
$A$-modules for $1 \le i \le p$.
\end{itemize}
The set of equivalence classes of such deformations is denoted
$\defm(R)$, and this defines the covariant functor $\defm$. Notice
that the fact that
    \[ M_R \cong (M_i \otimes_k R_{ij}) \]
as right $R$-modules replaces the flatness condition in classical
deformation theory. If $p=1$ and $R$ is commutative, the above
condition is of course equivalent to the flatness condition, so
the noncommutative deformation functor generalizes the classical
one.

In section \ref{s:mres}, we look at noncommutative deformations
from the point of view of resolutions. Let $R$ be any object of
$\acat p$. An M-free module over $R$ is a left $A \otimes_k
R^{\text{op}}$-module $F$ of the form
    \[ F = (L_i \otimes R_{ij}), \]
where $L_1, \dots, L_p$ are free left $A$-modules. M-free
complexes and M-free resolutions are defined similarly. Let us fix
a free resolution of $M_i$ the form
    \[ 0 \gets M_i \gets L_{0,i} \gets \dots \gets L_{m,i} \gets
    \dots \]
for $1 \le i \le p$. We prove that there is a bijective
correspondence between deformations of $\mfam$ to $R$ and
complexes of M-free modules over $R$ of the form
    \[ (L_{0,i} \otimes_k R_{ij}) \gets \dots \gets (L_{m,i}
    \otimes_k R_{ij}) \gets \dots \]
In fact, each such complex of M-free modules is an M-free
resolution of the corresponding deformation $M_R$ of $\mfam$ to
$R$.

In section \ref{s:hullf}, we recall some general facts about
pointed functors and their representability. In section
\ref{s:hulldef}, we consider the special case of the
noncommutative deformation functor $\defm$. From this point in the
text, we assume that the family $\mfam$ satisfy the finiteness
condition
\begin{equation}
\tag{\bfseries FC} \label{e:fc} \dim_k \ext^n_A(M_i,M_j) \text{ is
finite for } 1 \le i,j \le p, \; n=1,2.
\end{equation}
When this condition holds, we define $\infdef, \obstrdef$ to be
the formal matrix rings (in the sense of section \ref{s:cat})
given by the families of $k$-vector spaces $V_{ij} =
\ext^n_A(M_j,M_i)^*$ for $n = 1,2$. Assuming condition
(\ref{e:fc}), we show the following theorem of Laudal, which
generalizes the corresponding theorem for the classical
deformation functor:

\begin{thm}
There exists an obstruction morphism $o: \obstrdef \to \infdef$,
such that $H = \infdef \hat{\otimes}_{\obstrdef} k^p$ is a
pro-representable hull for the noncommutative deformation functor
$\defm: \acat p \to \sets$.
\end{thm}

In the rest of the paper, we show how to construct the hull $H$
explicitly, which can be accomplished by using matric Massey
products. In section \ref{s:mmp}, we introduce the immediately
defined matric Massey products. In section \ref{s:ch}, we define
the matric Massey products in general, and show that the hull $H$
of the noncommutative deformation functor $\defm$ is determined by
the vector spaces $\ext^n_A(M_i,M_j)$ for $n=1,2$ and $1 \le i,j
\le p$ and their matric Massey products. We also describe a
general method for calculating the hull $H$ in concrete terms.

In appendix \ref{s:app}, we describe the Yoneda and Hochschild
representations of the cohomology groups $\ext^n_A(M_i,M_j)$. In
this paper, we have chosen to express the matric Massey products
using the Yoneda representation and M-free resolutions. It is also
possible to express the matric Massey products using the
Hochschild representation, see for instance Laudal \cite{lau02}.

\section{Categories of pointed algebras}
\label{s:cat}

Let $p$ be a fixed natural number, and consider the ring $k^p$. This
ring has a natural $k$-algebra structure given by the map $\alpha
\mapsto (\alpha, \dots, \alpha)$ for $\alpha \in k$. Let $\pr_i: k^p
\to k^p$ be the $i$'th projection, and consider the ideal $k_i =
\pr_i(k^p) \subseteq k^p$ as a $k^p$-module for $1 \le i \le p$.
Clearly, $k^p$ is an Artinian $k$-algebra and $\{k_1, \dots, k_p \}$
is the full set of isomorphism classes of simple $k^p$-modules, each
of them of dimension $1$ over $k$. This simple example will serve as
a model for the $p$-pointed algebras that we shall consider in this
section.

A \emph{$p$-pointed $k$-algebra} is a triple $(R,f,g)$, where $R$ is
an associative ring and $f: k^p \to R, \; g: R \to k^p$ are ring
homomorphisms such that $g\circ f = \id$. A morphism $u: (R,f,g) \to
(R',f',g')$ of $p$-pointed $k$-algebras is a ring homomorphism $u: R
\to R'$ such that the natural diagrams commute (that is, such that $u
\circ f = f'$ and $g' \circ u = g$). We shall denote the category of
$p$-pointed $k$-algebras by $\Acat p$. Notice that if $(R,f,g)$ is an
object of $\Acat p$, then $f$ is injective and $g$ is surjective, and
we shall identify $k^p$ with its image in $R$. We often write $R$ for
the object $(R,f,g)$ to simplify notation.

Let $(R,f,g)$ be an object in $\Acat p$. We define the \emph{radical}
of $R$ to be $I(R) = \ker (g)$, which is an ideal in $R$. Furthermore,
we denote by $J(R)$ the \emph{Jacobson radical} of $R$
\[ J(R) = \{ x \in R: xM=0 \text{ for all simple left $R$-modules $M$}
\}, \]
which is also an ideal in $R$. We shall write $I,J$ for the radicals
$I(R),J(R)$ when there is no danger of confusion. Notice that the
Jacobson radical $J$ depends only on the ring $R$, while the radical
$I$ depends on the structural morphism $g$ as well.

For all objects $R$ in $\Acat p$, we have an inclusion $J(R)\subseteq
I(R)$: We have $J(k^p)=0$ since $k^p$ is semi-simple, and $g(J(R))
\subseteq J(k^p)=0$ since $g: R \to k^p$ is a surjection. In general,
we know that $R$ and $R/J(R)$ have the same simple left modules. So
if we consider $k_i$ as a left $R$-module via the morphism $g: R \to
k^p$ for $1 \le i \le p$, we see that $\{ k_1, \dots, k_p \}$ is
contained in the set of isomorphism classes of simple left $R$-modules,
and the equality $J(R)=I(R)$ holds if and only if $\{ k_1, \dots, k_p
\}$ is the full set of isomorphism classes of simple left $R$-modules.
Equivalently, the equality $I(R)=J(R)$ holds if and only if there are
exactly $p$ isomorphism classes of simple left $R$-modules.

It is therefore clear that the equality $I(R)=J(R)$ does not hold in
general: It  is easy to find examples where $R$ has `too many' simple
modules. For instance, consider $R=k[x]/(x-x^2)$ with the natural
$k$-algebra structure $f: k \to R$ and let $g: R \to k$ be given by $x
\mapsto 0$. Then $R$ is an object of $\Acat 1$, but $J(R) \neq I(R)$
because $R$ has two non-isomorphic simple left $R$-modules (given by $x
\mapsto 0$ and $x \mapsto 1$).

Let $e_i$ be the idempotent $(0,0, \dots , 1, \dots, 0) \in k^p$ for $1
\le i \le p$. Notice that $e_i e_j = 0$ if $i\neq j$, and that
$e_1 + \dots + e_p = 1$. For any object $R$ in $\Acat p$, we identify
$\{ e_1, \dots, e_p \}$ with idempotents in $R$ via the inclusion $k^p
\to R$. Denote by $R_{ij}$ the $k$-linear sub-space $e_i R e_j \subseteq
R$. We immediately see, using the properties of the idempotents, that
the following relations hold for $1 \le i,j,l,m \le p$:
\begin{enumerate}
\item $R_{ij} R_{lm} \subseteq \delta_{jl} R_{im}$,
\item $R_{ij} \cap R_{lm} = 0$ if $(i,j) \neq (l,m)$,
\item $\sum R_{ij} = R$.
\end{enumerate}
In particular, we have that $R= \oplus R_{ij}$, so every element $r \in
R$ may be written in matrix form $r=(r_{ij})$ with $r_{ij} \in R_{ij}$
for $1 \le i,j \le p$. Furthermore, elements of $R$ multiply as matrices
when we write them in this form. It is therefore reasonable to call an
object $R$ in $\Acat p$ a matrix ring, and to write it $R=(R_{ij})$.
Notice that $R_{ii}$ is an associative ring (with identity $e_i$), and
that $R_{ij}$ is a (unitary) $R_{ii} - R_{jj}$ bimodule for $1\le i,j
\le p$. For any ideal $K \subseteq R$, we see that $e_i K e_j = K \cap
R_{ij}$, and we shall denote this $k$-linear subspace $K_{ij}$ for $1\le
i,j \le p$. Since $K=\osum K_{ij}$, we write $K=(K_{ij})$.

Let $R$ be an object of $\Acat p$, so $R=(R_{ij})$ is a matrix ring in the
above sense. The following standard result gives useful information on
when $R$ is an Artinian or Noetherian ring:

\begin{prop}
\label{p:mcc}
Let $R=(R_{ij})$ be an object in $\Acat p$. Then $R$ is Noetherian
(Artinian) if and only if the following conditions hold:
\begin{thmlist}
\item $R_{ii}$ is Noetherian (Artinian) for $1\le i \le p$,
\item $R_{ij}$ is a Noetherian (Artinian) left $R_{ii}$-module and a
Noetherian (Artinian) right $R_{jj}$-module for $1\le i \neq j \le p$.
\end{thmlist}
\end{prop}

We recall that a finitely generated, associative $k$-algebra is not
necessarily Noetherian. That is, Hilbert's basis theorem does not hold
for associative rings. For a counter-example, let $R=k\{ x_1,\dots,
x_n \}$ be the free associative $k$-algebra on $n$ generators. It is
well-known that $R$ is Noetherian only if $n=1$. However, we know from
the Hopkins-Levitzki theorem that an associative Artinian ring is
Noetherian.

A $k$-algebra $R$ of finite dimension as vector space over $k$ is
Artinian. This is clear, since every one-sided ideal is a vector space
over $k$ of finite dimension. We have a converse statement under
the following conditions:

\begin{lem}
\label{l:fd}
Let $R$ be an object of $\Acat p$. If $R$ is Artinian and $I(R)$ is
nilpotent, then $R$ has finite dimension as a vector space over $k$.
\end{lem}
\begin{proof} \noindent
We write $I = I(R)$. Since $R$ is Artinian and therefore Noetherian,
$I^m$ is finitely generated as a left $R$-module for all $m$.
Consequently, $I^m/I^{m+1}$ is a finitely generated $R/I$-module for
all $m$, and hence has finite $k$-dimension. But $I^n=0$ for some $n$,
so $I^m$ has finite $k$-dimension for all $m \ge 0$. In particular,
$R$ has finite dimension as a vector space over $k$.
\end{proof}

We define the \emph{category} $\acat p$ to be the full sub-category
of $\Acat p$ consisting of objects $R$ in $\Acat p$ such that $R$ is
Artinian and $I(R)=J(R)$. The condition $I(R)=J(R)$ might
equivalently be replaced by the condition that $I(R)$ is a nilpotent
ideal, since the Jacobson radical is the largest nilpotent ideal in
an Artinian ring. So by lemma \ref{l:fd}, all objects $R$ in $\acat
p$ have finite $k$-dimension. Since $R$ is Artinian, the condition
that $I(R)$ is nilpotent is also equivalent to $\cap \; I(R)^n = 0$.
Finally, there is a geometric interpretation of the condition
$I(R)=J(R)$: By the comment earlier in this section, $I(R)=J(R)$ if
and only if $\{ k_1, \dots, k_p \}$ is the full set of isomorphism
classes of simple left $R$-modules (or equivalently, that the number
of such isomorphism classes is exactly $p$).

\begin{lem}
Let $R$ be an associative ring. Then there exists morphisms $f: k^p
\to R$ and $g: R \to k^p$ making $(R,f,g)$ an object of $\acat p$ if
and only if $R$ is an Artinian $k$-algebra with exactly $p$
isomorphism classes of simple left $R$-modules, each of them of
dimension $1$ over $k$.
\end{lem}
\begin{proof}
One implication follows from the comments above. For the other, assume
that $R$ is Artinian with the prescribed isomorphism classes of simple
left $R$-modules. This defines a morphism $g: R \to k^p$. Clearly, $I=
\ker(g) = J(R)$ by the comments above. So it is enough to lift the
idempotents $\{ e_1, \dots, e_p \}$ of $k^p$ to idempotents $\{ r_1,
\dots, r_p \}$ in $R$ such that $r_1 + \dots + r_p = 1$ and $r_i r_j =
0$ when $i \neq j$. But $R$ is Artinian and therefore $I = J(R)$ is
nilpotent, so this is clearly possible.
\end{proof}

Let $R$ be an object in $\Acat p$ with radical $I=I(R)$. Then the $I$-adic
filtration defines a topology on $R$ compatible with the ring operations,
and we shall always consider $R$ a topological ring in this way. We say
that the topology on $R$ is Hausdorff (or separated) if and only if
$\cap I^n =0$.

For all objects $R$ in $\Acat p$, there is an $I$-adic completion $\hat{R}$
of $R$ and a canonical morphism $R \to \hat{R}$ in $\Acat p$. The $I$-adic
completion $\hat{R}$ is defined by the projective limit
\[ \hat{R} = \lim_{\leftarrow} R/I^n, \]
and the morphism $R \to \hat{R}$ is the natural one induced by
this projective limit. Notice that the kernel of this morphism is
$\cap I^n$. We say that $R$ is \emph{complete} (or separated
complete) if the natural morphism $R \to \hat{R}$ is an
isomorphism in $\Acat p$. In particular, this implies that the
morphism is injective, so $R$ is Hausdorff (or separated). This
gives a new characterization of the category $\acat p$:

\begin{lem}
The category $\acat p$ is the full sub-category of $\Acat p$
consisting of objects such that $R$ is Artinian and $I$-adic
complete.
\end{lem}

We define the \emph{pro-category} $\cacat p$ of $\acat p$ to be
the full sub-category of $\Acat p$ consisting of objects such that
$R$ is complete and $R/I(R)^n$ belongs to $\acat p$ for all $n\ge
1$. It is clear that we have an inclusion of (full) sub-catgories
$\acat p \subseteq \cacat p$.

Let $R$ be an object in $\cacat p$ with radical $I=I(R)$. To fix notation,
we write $\gr_n(R) = I^n / I^{n+1}$ for $n \ge 0$ (with $I^0 = R$). We also
write $\gr R = \osum \gr_n(R)$, this is the graded ring associated
to the $I$-adic filtration of $R$. The \emph{tangent space} of $R$ is
defined to be the $k$-linear space dual to $\gr_1(R)$,
\[ \tg R = \hmm_k(I/I^2,k) = (I/I^2)^*, \]
which is clearly of finite dimension over $k$. In particular, we have
$(\tg R)^* \cong I/I^2$.

Let $u: R \to S$ be a morphism in $\cacat p$. As usual, we consider $R$
and $S$ with the $I$-adic filtrations, where $I$ is $I(R)$ and $I(S)$
respectively. Since $u$ preserves these filtrations, it induces a
morphism of graded rings $\gr(u): \gr R \to \gr S$. This morphism is
homogeneous of degree $0$, so $u$ also induces morphisms of $k$-vector
spaces $\gr_n(u): \gr_n(R) \to \gr_n(S)$ for all $n \ge 0$. In particular,
we have a morphism of $k$-vector spaces $\gr_1(u): \gr_1(R) \to \gr_1(S)$,
and a dual morphism $\tg u: \tg S \to \tg R$.

\begin{prop}
\label{p:surj} Let $u:R \to S$ be a morphism in $\cacat p$. Then
$u$ is a surjection if and only if $\gr_1(u)$ is a surjection.
Furthermore, u is injective if $\gr (u)$ is injective.
\end{prop}
\begin{proof} \noindent
If $u$ is surjective, then clearly $\gr_1(u)$ is also surjective.
To prove the other implication, let us consider the map $\gr(u):
\gr(R) \to \gr(S)$. Since $\gr S$ is generated by the elements in
$\gr_1 S$ as an algebra, it follows that if $\gr_1 (u)$ is
surjective, then $\gr (u)$ is also surjective. From Bourbaki
\cite{bou85}, chapter III, \S 2, no. 8, corollary 1 and 2, we have
that $u$ is surjective (injective) if $\gr (u)$ is surjective
(injective), and the result follows.
\end{proof}

Let $n$ be any natural number. We define the \emph{category} $\acat p(n)$
to be the full sub-category of $\acat p$ consisting of objects $R$ in
$\acat p$ such that $I(R)^n=0$. Notice that $\acat p(n) \subseteq \acat
p(n+1)$ for all $n\ge 1$. Furthermore, each object $R$ in $\acat p$
belongs to a sub-category $\acat p(n)$ for some integer $n$.

Let $u:R \to S$ be a morphism in $\acat p$, and denote by
$K=\ker(u)$ the kernel of $u$. We say that $u$ is a \emph{small
morphism} if we have $I(R) \cdot K = K \cdot I(R) = 0$. We prove
the following important fact about small surjections:

\begin{lem}
Let $u:R \to S$ be a surjection in $\acat p$. Then $u$ can be
factored into a finite number of small surjections.
\end{lem}
\begin{proof}
Let $I = I(R)$, then $I^n K = 0$ for some $n \ge 0$. Consider the
surjection $u_q: R/I^q K \to R/I^{q-1} K$ for $1 \le q \le n$.
Clearly $I(R/I^q K) \ker(u_q) = 0$ for all $q$. Moreover, $u_1
\circ \dots \circ u_n = u$ when $u_1: R/IK \to R/K$ is considered
as a morphism onto $S \cong R/K$. It is therefore enough to prove
the lemma for a surjection $u:R \to S$ with $IK=0$. In this
situation, $K I^n = 0$ for some $n \ge 0$. Now consider the
surjection $v_q: R/K I^q \to R/K I^{q-1}$ for $1 \le q \le n$.
Clearly, $v_q$ is a small surjection for all $q$. Moreover, $u =
v_1 \circ \dots \circ v_n$ when $v_1: R/KI \to R/K$ is considered
as a morphism onto $S \cong R/K$. It follows that $u$ can be
factorized in a finite number of small surjections in $\acat p$.
\end{proof}

We conclude this section with an important family of examples:
Let $V_{ij}$ be a finite dimensional $k$-vector space for
$1\le i,j \le p$, with $\dim_k V_{ij} = d_{ij}$. Let furthermore
$\{ r_{ij}(l): 1 \le l \le d_{ij} \}$ be a basis of $V_{ij}$ for
$1\le i,j \le p$ (or simply $\{ r_{ij} \}$ if $d_{ij}=1$). We
define the \emph{free matrix ring} $R=R(\{V_{ij}\})$ defined by
the vector spaces $V_{ij}$ in the following way: We say that a
monomial in $R$ of type $(i,j)$ and degree $n$ is an expression
of the form
\[ r_{i_0 i_1}(l_1) r_{i_1 i_2}(l_2) \dots r_{i_{n-1} i_n}(l_n) \]
with $i_0=i, i_n=j$. To these, we add the monomials $e_i$ for
$1\le i\le p$, which we consider to be of type $(i,i)$ and degree
$0$. We define $R$ to be the $k$-linear space generated by all
monomials in $R$, with the obvious multiplication: If $M$ is a
monomial of type $(i,j)$, and $M'$ is a monomial of type $(l,m)$,
then $M M' = 0$ if $j \neq l$, and $M M'$ is the monomial obtained
by juxtapositioning $M$ and $M'$ (possibly after having erased
unnecessary $e_i$'s) if $j=l$. We see that $(R,f,g)$ is an object
of the category $\Acat p$, where $f,g$ are the obvious maps $k^p
\to R \to k^p$. In fact, $R_{ij}$ is the $k$-linear subspace
generated by monomials in $R$ of type $(i,j)$, and the ideal
$I=I(R)$ is the $k$-linear subspace generated by all monomials of
positive degree.

We denote by $\hat{R}=\hat{R}(\{V_{ij}\})$ the completion of $R =
R(\{V_{ij}\})$, and call this the \emph{formal matrix ring}
defined by the vector spaces $V_{ij}$. Explicitly, every element
in $\hat{R}_{ij}$ is an infinite $k$-linear sum of monomials in
$R$ of type $(i,j)$. Let $I=I(R)$, then we have that $R_n = R/I^n
\cong \hat{R}/I(\hat{R})^n$ belongs to $\acat p$ for $n\ge 1$:
Clearly, $R_n$ has finite dimension as $k$-vector space, so $R_n$
is Artinian, and $I(R_n)=I/I^n$, so the radical is nilpotent.
Since $\hat{R}$ clearly is complete, it follows that $\hat{R}$
belongs to $\cacat p$.

Notice that neither the free matrix ring $R$ nor the formal
matrix ring $\hat{R}$ is Noetherian in general. For a
counter-example, it is enough to consider the case when $p=1$ and
$d_{11}=2$, or the case when $p=2$ and $d_{11}=d_{12}=d_{21}=1,
\; d_{22}=0$. In the first case, $R \cong k\{ x,y \}$, which we
know is not Noetherian. In the second case, we have that $R_{11}=
k\{r_{11}, r_{12} r_{21}\} \cong k\{x,y\}$, which again is not
Noetherian. So by proposition \ref{p:mcc}, $R$ is not Noetherian
in this case either. A similar argument shows that $\hat R$ is
not Noetherian in any of the two cases.

\section{Noncommutative deformations of modules}
\label{s:ncdefm}

We recall that $k$ is an algebraically closed (commutative) field,
$A$ is an associative $k$-algebra, and $\mfam =\{ M_1, \dots, M_p \}$
is a finite family of left $A$-modules. In this section, we shall
define the noncommutative deformation functor
\[ \defm: \acat p \to \sets \]
describing the simultaneous formal deformations of the family
$\mfam$.

Let $R$ be an object of $\acat p$. A lifting of the family $\mfam$
of left $A$-modules to $R$ is a left $A\otimes_k R^{\text{op}}$-module
$M_R$, together with isomorphisms $\eta_i: M_R \otimes_R k_i \to M_i$
of left $A$-modules for $1\le i \le p$, such that $M_R \cong (M_i
\otimes_k R_{ij})$ as right $R$-modules. We remark that a left
$A\otimes_k R^{\text{op}}$-module is the same as an $A$-$R$ bimodule
such that the left and right $k$-vector space structures coincide.
Furthermore, the notation $(M_i \otimes_k R_{ij})$ refers to the
$k$-vector space
\[ (M_i \otimes_k R_{ij}) = \osum_{i,j} (M_i \otimes_k
R_{ij}) \]
with the natural right $R$-module structure coming from the
multiplication in $R$. The condition that $M_R \cong (M_i \otimes_k
R_{ij})$ as right $R$-modules generalizes the flatness condition in
commutative deformation theory.

Let $M'_R, M''_R$ be two liftings of $\mfam$ to $R$. We say that
these two liftings are equivalent if there exists an isomorphism
$\tau: M'_R \to M''_R$ of left $A\otimes_k R^{\text{op}}$-modules
such that the natural diagrams commute (that is, such that $\eta''_i
\circ (\tau \otimes_R k_i) = \eta'_i$ for $1\le i \le p$). We let
$\defm (R)$ denote the set of equivalence classes of liftings of
$\mfam$ to $R$, and we refer to these equivalence classes as
\emph{deformations} of $\mfam$ to $R$. We shall often denote a
deformation represented by $(M_R, \eta_i)$ by $M_R$ to simplify
notation.

Let $u: R \to S$ be a morphism in $\acat p$, and let $M_R$ be
a lifting of $\mfam$ to $R$, representing an element in $\defm
(R)$. We define $M_S = M_R \otimes_R S$, which has a natural
structure as a left $A\otimes_k S^{\text{op}}$-module. Since
$u$ is a morphism in $\acat p$, we have natural isomorphisms
of left $A$-modules
\[ (M_R \otimes_R S) \otimes_S k_i \cong M_R \otimes_R k_i, \]
inducing isomorphisms of left $A$-modules $\rho_i: M_S \otimes_S k_i
\to M_i$ via $\eta_i$ for $1\le i \le p$. A straight-forward
calculation shows that $M_S$ together with the isomorphisms $\rho_i$
for $1\le i \le p$ constitutes a lifting of $\mfam$ to $S$, and
furthermore that the equivalence class of this lifting is independent
upon the representative of the equivalence class of $M_R$. Hence, we
obtain a map $\defm(u): \defm(R) \to
\defm(S)$, and we see that $\defm: \acat p \to \sets$ is a
covariant functor.

Let $R=(R_{ij})$ be an object in $\acat p$. We shall describe how
one, in principle, could attempt to calculate $\defm(R)$ explicitly:
We may assume that every element of $\defm(R)$ is represented by a
lifting $M_R$, such that $M_R=(M_i \otimes_k R_{ij})$ considered as
a right $R$-module. In order to describe this lifting completely,
it is enough to describe the left action of $A$ on $M_R$.
Furthermore, it is enough to describe this action on elements
of the form $m_i \otimes e_i$ with $m_i \in M_i$, since we have
\[ a ( m_i \otimes r_{ij}) = (a(m_i \otimes e_i)) r_{ij} \]
for all $a \in A, \; m_i \in M_i, \; r_{ij} \in R_{ij}$. For a
fixed $a \in A, \; m_i \in M_i$, assume that $a(m_i \otimes e_i) =
\sum (m'_j \otimes r'_{jl})$ with $m'_j \in M_j, \; r'_{jl} \in
R_{jl}$. Then multiplication by $e_i$ on the right gives the
equality
\[ a(m_i \otimes e_i) = \sum_{j} (m'_j \otimes r'_{ji}), \]
and the isomorphism $\eta_i$ gives a further restriction on the
left action of $A$, expressed by the formula
\begin{equation}
\label{f:fdef}
a(m_i \otimes e_i) = (a m_i) \otimes e_i + \sum_j m'_j \otimes
r'_{ji},
\end{equation}
where $a \in A, \; m_i \in M_i, \; m'_j \in M_j, \; r'_{ji} \in
I(R)_{ji}$. Consequently, the set $\defm(R)$ consists of all
possible choices of left $A$-actions on elements of the form $m_i
\otimes e_i$, fulfilling condition (\ref{f:fdef}) and the
associativity condition, up to equivalence.

Let $R$ be any object in $\acat p$. Then the formula $a(m_i \otimes
e_i)=(am_i) \otimes e_i$ for $a \in A, \; m_i \in M_i$ defines a left
$A$-module structure on $(M_i \otimes R_{ij})$ compatible with the
right $R$-module structure. Hence, there exists a trivial lifting
$M_R$ to $R$ for all objects $R$ in $\acat p$, and $\defm(R)$ is
non-empty. Notice that in the case $R=k^p$, we have $I=I(R)=0$, so
this trivial lifting is the only one possible. Consequently, we have
$\defm(k^p)= \{ * \}$, where $*$ denotes the equivalence class of the
trivial lifting.

Let $u:R \to S$ be a morphism in $\acat p$, and let $M_S \in \defm(S)$
be a given deformation. We say that a deformation $M_R \in \defm(R)$
is a lifting of $M_S$ or is lying over $M_S$ if $\defm(u)(M_R)=M_S$.
Given any object $R$ in $\acat p$ and a deformation $M_R \in \defm(R)$,
we see that $M_R$ is a lifting of the trivial deformation $*$ in
$\defm(k^p)$ in the above sense via the structural morphism $g:R \to
k^p$. Hence, our notation is consistent.

For another example, consider the test algebras $\test{\alpha}{\beta}$
for $1 \le \alpha, \beta \le p$, constructed in the following way: Let
$R$ be the free matrix algebra defined by the $k$-vector spaces $V_{ij}$
with dimensions $d_{\alpha,\beta}=1$ and $d_{ij}=0$ when $(i,j) \ne
(\alpha,\beta)$. We define $\test{\alpha}{\beta}=R/I(R)^2$, which is an
object in $\acat p(2)$ by construction. We know that any lifting of
$\mfam$ to $\test{\alpha}{\beta}$ is defined by a left $A$-action
\[ a(m_{\beta} \otimes e_{\beta}) = (am_{\beta}) \otimes e_{\beta} +
\psi(a)(m_{\beta}) \otimes \varepsilon_{\alpha,\beta} \]
for all $a \in A, \; m_{\beta} \in M_{\beta}$, where $\psi: A \times
M_{\beta} \to M_{\alpha}$ is a $k$-bilinear map and $\varepsilon_{\alpha,
\beta}$ is the class of $r_{\alpha,\beta}$. Clearly, we must have
$a(m_i \otimes e_i) = (am_i) \otimes e_i$ for all $a \in A, m_i \in M_i$
when $i \neq \beta$. Moreover, $\psi$ defines an associative $A$-module
structure if and only if $\psi \in \der_k(A,\hmm_k(M_{\beta},
M_{\alpha}))$. In this case, we shall denote the corresponding lifting
by $\mfam(\psi) \in \defm(\test{\alpha}{\beta})$. Given two derivations
$\psi, \psi'$, we see that $\mfam(\psi)$ and $\mfam(\psi')$ are
equivalent liftings if and only if there is a $\phi \in \hmm_k(M_{\beta},
M_{\alpha})$ such that
\[ (\psi-\psi')(a)(m_{\beta}) =  a\phi(m_{\beta}) - \phi(am_{\beta}) \]
for all $a \in A, \; m_{\beta} \in M_{\beta}$.

\begin{lem}
\label{l:testdefm}
There is a bijective correspondence $\defm(\test{\alpha}{\beta}) \cong
\ext^1_A(M_{\beta}, M_{\alpha})$ for $1 \le \alpha, \beta \le p$.
\end{lem}
\begin{proof}
From the definition of Hochschild cohomology (see appendix
\ref{a:hh}), we see that $\psi \mapsto \mfam(\psi)$ induces a
bijective correspondence between
$\hh^1(A,\hmm_k(M_{\beta},M_{\alpha}))$ and
$\defm(\test{\alpha}{\beta})$. Moreover,
$\hh^1(A,\hmm_k(M_{\beta}, M_{\alpha})) \cong
\ext^1_A(M_{\beta},M_{\alpha})$ by proposition \ref{p:hhext}.
\end{proof}

\section{M-free resolutions and noncommutative deformations}
\label{s:mres}

We recall that $k$ is an algebraically closed (commutative) field,
$A$ is an associative $k$-algebra, and $\mfam =\{ M_1, \dots, M_p
\}$ is a finite family of left $A$-modules. In this section, we
shall define M-free resolutions and relate them to noncommutative
deformations of modules. In particular, we shall show that M-free
resolutions are useful computational tools in order to study the
deformation functor $\defm$.

Let $R$ be any object of $\acat p$. An \emph{M-free module} over $R$
is a left $A \otimes_k R^{\text{op}}$-module $F$ of the form
    \[ F = (L_i \otimes_k R_{ij}), \]
where $L_1, \dots, L_p$ are free left $A$-modules, and the left
$A$-module structure on $F$ is the trivial one. In other words, $F$
is the trivial lifting of a family $\{ L_1, \dots, L_p \}$ of free
left $A$-modules to $R$.

Although an M-free module over $R$ is not free considered as a
left $A \otimes_k R^{\text{op}}$-module, it behaves as a free
module when interpreted as a module of matrices in the correct
way:

\begin{lem}
\label{l:mflift} Let $u: R \to S$ be a surjection in $\acat p$,
and consider a left $A \otimes_k R^{\text{op}}$-module $M_R = (M_i
\otimes_k R_{ij})$ and a left $A \otimes_k S^{\text{op}}$-module
$M_S = (M_i \otimes_k S_{ij})$ such that the natural map $v: M_R
\to M_S$ induced by $u$ is left $A$-linear. If $F^S$ is any M-free
module over $S$ given by the free left $A$-modules $L_1, \dots,
L_p$ and $f_S: F^S \to M_S$ is any left $A \otimes_k
S^{\text{op}}$-linear map, then there exists a left $A \otimes_k
R^{\text{op}}$-linear map $f_R: F^R \to M_R$ making the diagram
\[
\xymatrix{
M_R \ar[d]_{v} & F^R \ar[d]^{(\id \otimes u)} \ar[l]_{f_R} \\
M_S & F^S \ar[l]^{f_S} \\
}
\]
commutative, where $F^R$ is the M-free module over $R$ given by
the free left $A$-modules $L_1, \dots, L_p$.
\end{lem}
\begin{proof}
Clearly, the map $f_S$ is determined by its values on $L_i \otimes e_i$,
and therefore by the corresponding left $A$-linear maps $L_j \to \oplus
(M_i \otimes_k S_{ij})$. Since each left $A$-module $L_j$ is projective,
we can lift these maps to left $A$-linear maps $L_j \to \oplus (M_i
\otimes_k R_{ij})$, and these maps determine $f_R$.
\end{proof}

Let $R$ be any object of $\acat p$, and let $M_R = (M_i \otimes_k R_{ij})
\in \defm(R)$ be a lifting of $\mfam$ to $R$. An \emph{M-free resolution}
of $M_R$ is an exact sequence of left $A \otimes_k R^{\text{op}}$-linear
maps
\[ 0 \gets M_R \gets F^R_0 \gets F^R_1 \gets \dots \gets F^R_m \gets
\dots \]
where $F^R_m$ is an M-free module over $R$ for $m \ge 0$. So we have
$F^R_m = (L_{m,i} \otimes_k R_{ij})$ where $L_{m,i}$ are free left
$A$-modules for $1 \le i \le p, \; m \ge 0$. We shall denote the
differentials by $d^R_m: F^R_{m+1} \to F^R_m$ for $m \ge 0$.

We fix a $k$-linear basis $\{ r_{ij}(l): 1 \le l \le \dim_k R_{ij} \}$
of $R_{ij}$ for $1 \le i,j \le p$ such that $e_i$ is contained in the
basis of $R_{ii}$ for $1 \le i \le p$. Consider the differential $d^R_m$
in the M-free resolution of $M_R$ above. Clearly, we can write $d^R_m$
uniquely in the form
\begin{equation}
\label{e:diffrep}
d^R_m = \sum_{i,j,l} \; \alpha(r_{ij}(l))_m \otimes r_{ij}(l)
\end{equation}
for all $m \ge 0$, where $\alpha(r_{ij}(l))_m: L_{m+1,j} \to L_{m,i}$ is
a homomorphism of left $A$-modules for $1 \le i,j \le p, \; 1 \le l \le
\dim_k R_{ij}$. In particular, the M-free resolution of $M_R$ defines a
family of $1$-cochains $\alpha(r_{ij}(l)) \in \hmm^1(L_{\cpd j},
L_{\cpd i})$, indexed by a $k$-linear basis for $R$.

From now on, we fix a free resolution $(L_{\cpd i},d_{\cpd i})$ of
$M_i$ considered as left $A$-module for $1 \le i \le p$. These
free resolutions correspond to an M-free resolution
$(F_{\cpd},d_{\cpd})$ of the trivial deformation $(M_i \otimes_k
(k^p)_{ij}) \in \defm(k^p)$. In fact, the M-free resolution
$(F_{\cpd},d_{\cpd})$ is given by $F_m = (L_{m,i} \otimes_k
(k^p)_{ij})$ and $d_m = \sum d_{m,i} \otimes e_i$ for $m \ge 0$.
We have therefore fixed an M-free resolution $(F_{\cpd},d_{\cpd})$
of the trivial lifting of $\mfam$ to $k^p$.

Let $R$ be any object of $\acat p$. We say that a complex $(F^R_{\cpd},
d^R_{\cpd})$ of M-free modules $F^R_m = (L_{m,i} \otimes_k R_{ij})$ over
$R$ is a \emph{lifting of the complex $(F_{\cpd},d_{\cpd})$} if the
following diagram commutes
\[
\xymatrix{
F^R_0 \ar[d]_{v_0} & F^R_1 \ar[l]_{d^R_0} \ar[d]_{v_1} & F^R_2
\ar[l]_{d^R_1} \ar[d]_{v_2} & \dots \ar[l] \\
F_0 & F_1 \ar[l]^{d_0} & F_2 \ar[l]^{d_1} & \dots \ar[l] \\
}
\]
where $v_m: F^R_m \to F_m$ are the natural maps induced by $R \to k^p$.

\begin{lem}
\label{l:mfcp}
Let $R$ be any object of $\acat p$, and let $(F^R_{\cpd},d^R_{\cpd})$
be a lifting of the complex $(F_{\cpd},d_{\cpd})$. Then we have:
\begin{enumerate}
\item $H^m(F^R_{\cpd},d^R_{\cpd}) = 0$ for all $m \ge 0$,
\item $H^0(F^R_{\cpd},d^R_{\cpd})$ is a lifting of the family $\mfam$
to $R$.
\end{enumerate}
\end{lem}
\begin{proof}
Clearly, the lemma holds for $R=k^p$. We shall consider a small
surjection $u:R \to S$ in $\acat p$ and liftings of complexes
$(F^U_{\cpd},d^U_{\cpd})$ of $(F_{\cpd},d_{\cpd})$ to $U$ for $U=R,S$
such that the following diagram commutes:
\[
\xymatrix{
F^R_0 \ar[d]_{v_0} & F^R_1 \ar[l]_{d^R_0} \ar[d]_{v_1} & F^R_2
\ar[l]_{d^R_1} \ar[d]_{v_2} & \dots \ar[l] \\
F^S_0 & F^S_1 \ar[l]^{d^S_0} & F^S_2 \ar[l]^{d^S_1} & \dots \ar[l] \\
}
\]
In this situation, we shall prove that if the conclusion of the
lemma holds for $S$, it holds for $R$ as well. This is clearly
enough to prove the lemma.

Let $K=\ker(u)$, then we clearly have $\ker(v_m) = (F_{m,i} \otimes_k
K_{ij})$ with the trivial left $A$-action for all $m \ge 0$. We denote
this kernel by $F^K_m$, then $(F^K_{\cpd},d^K_{\cpd})$ is a complex of
left $A \otimes_k R^{\text{op}}$-modules, where $d^K_{\cpd}$ is the
restriction of $d^R_{\cpd}$. Moreover, it is clear that $v_m$ is
surjective for $m \ge 0$. Define $M_U = H^0(F^U_{\cpd},d^U_{\cpd})$ for
$U=R,S$, let $v: M_R \to M_S$ be the induced map, and denote the kernel
by $M_K = \ker(v)$. Then clearly $v$ is surjective, and we have the
following commutative diagram of complexes:
\[
\xymatrix{
& 0 \ar[d] & 0 \ar[d] & 0 \ar[d] & 0 \ar[d] & \\
0 & M_K \ar[l] \ar[d]_i & F^K_0 \ar[l]_{\rho^K} \ar[d]_{i_0} & F^K_1
\ar[l]_{d^K_0} \ar[d]_{i_1} & F^K_2 \ar[l]_{d^K_1} \ar[d]_{i_2} & \dots
\ar[l] \\
0 & M_R \ar[l] \ar[d]_v & F^R_0 \ar[l]_{\rho^R} \ar[d]_{v_0} & F^R_1
\ar[l]_{d^R_0} \ar[d]_{v_1} & F^R_2 \ar[l]_{d^R_1} \ar[d]_{v_2} & \dots
\ar[l] \\
0 & M_S \ar[d] \ar[l] & F^S_0 \ar[d] \ar[l]^{\rho^S} & F^S_1 \ar[d]
\ar[l]^{d^S_0} & F^S_2 \ar[d] \ar[l]^{d^S_1} & \dots \ar[l] \\
& 0 & 0 & 0 & 0 & \\
}
\]
Clearly all columns are exact, so the diagram gives a short exact
sequence of complexes. By assumption, the bottom row is exact and
$M_S = (M_i \otimes_k S_{ij})$ is a lifting of $\mfam$ to $S$. Let
us first show that $H^m(F^K_{\cpd},d^K_{\cpd})=0$ for $m \ge 1$: This
follows since the complex is a lifting of $(F_{\cpd},d_{\cpd})$ and
because $I(R) K=0$ (since $u:R \to S$ is small). The long exact
sequence of cohomologies of the complexes above now implies that
$H^m(F^R_{\cpd},d^R_{\cpd})=0$ for all $m \ge 1$ and that we have a
short exact sequence
\[ 0 \to H^0(F^K_{\cpd},d^K_{\cpd}) \to M_R \to M_S = (M_i \otimes_k
S_{ij}) \to 0, \]
of left $A$-modules, so in particular $M_K \cong H^0(F^K_{\cpd},
d^K_{\cpd})$. But since $I(R) K = 0$, it follows that $H^0(F^K_{\cpd},
d^K_{\cpd}) \cong (H^0(L_{\cpd,i},d_{\cpd,i}) \otimes_k K_{ij}) =
(M_i \otimes_k K_{ij})$ with the trivial left $A$-module structure. It
follows that $M_R \cong (M_i \otimes_k R_{ij})$ considered as a
$k$-vector space, and therefore $M_R$ is a lifting of $\mfam$ to $R$.
\end{proof}

\begin{lem}
\label{l:mfres}
Let $R$ be any object of $\acat p$, and let $M_R$ be a lifting of
$\mfam$ to $R$. Then there exists an M-free resolution of $M_R$ which
lifts the complex $(F_{\cpd},d_{\cpd})$ to $R$.
\end{lem}
\begin{proof}
Clearly, the lemma holds for $R=k^p$. We shall consider a small
surjection $u:R \to S$ in $\acat p$, deformations $M_U \in \defm(U)$
for $U=R,S$ such that $M_R$ lifts $M_S$ to $R$, and an M-free
resolution $(F^S_{\cpd},d^S_{\cpd})$ of $M_S$ which lifts the complex
$(F_{\cpd},d_{\cpd})$ to $S$. In this situation, we shall prove that
there exists an M-free resolution $(F^R_{\cpd},d^R_{\cpd})$ of $M_R$
compatible with the M-free resolution of $M_S$. This is clearly enough
to prove the lemma.

Let $F^R_m = (L_{m,i} \otimes_k R_{ij})$ for all $m \ge 0$. Moreover,
we write $F^K_m = (L_{m,i} \otimes_k K_{ij})$ for all $m \ge 0$, where
$K = \ker(u)$. To complete the proof, we have to find the differentials
$d^R_m$ for $m \ge 0$ and the augmentation map $\rho_R$: By lemma
\ref{l:mflift}, we can find a homomorphism $\rho_R: F^R_0 \to M_R$
lifting $\rho_S$. Denote by $\rho_K: F^K_0 \to M_K$ its restriction,
where $M_K = \ker(M_R \to M_S)$. Since $u$ is small, $\rho^K$ is
surjective, and this implies that the induced map $\ker(\rho^R) \to
\ker(\rho^S)$ is surjective. By lemma \ref{l:mflift}, we can find a
homomorphism $d^R_0: F^R_1 \to F^R_0$ lifting $d^S_0$ such that $\rho^R
d^R_0 = 0$. Let $d^K_0$ be the restriction of $d^R_0$, then clearly
$\ker(\rho^K)=\im(d^K_0)$ since $u$ is small. An easy induction argument
shows that we can construct a complex $(F^R_{\cpd},d^R_{\cpd})$ lifting
the complex $(F^S_{\cpd},d^S_{\cpd})$ in such a way that the restriction
$(F^K_{\cpd},d^K_{\cpd})$ is a resolution of $M_K$. By the proof of
lemma \ref{l:mfcp}, it follows that $H^m(F^R_{\cpd},d^R_{\cpd})=0$ for $m
\ge 1$ and that there is an exact sequence
\[ 0 \to M_K \to H^0(F^R_{\cpd},d^R_{\cpd}) \to M_S \to 0. \]
This implies that $M_R = H^0(F^R_{\cpd},d^R_{\cpd})$, and $(F^R_{\cpd},
d^R_{\cpd})$ is the required M-free resolution of $M_R$ compatible with
the given M-free resolution of $M_S$.
\end{proof}

\begin{prop}
\label{p:reslift} Let $u:R \to S$ be a surjection in $\acat p$,
and consider a deformation $M_S \in \defm(S)$ and any M-free
resolution $(F^S_{\cpd},d^S_{\cpd})$ of $M_S$ which lifts the
complex $(F_{\cpd},d_{\cpd})$ to $S$. There is a bijective
correspondence between the set of liftings
    \[ \{ M_R \in \defm(R): \defm(u)(M_R) = M_S \} \]
and the set of M-free complexes $(F^R_{\cpd},d^R_{\cpd})$ which
lift the resolution $(F^S_{\cpd},d^S_{\cpd})$ to $R$, up to
equivalence.
\end{prop}
\begin{proof}
For a small surjection, this follows from lemma \ref{l:mfcp} and lemma
\ref{l:mfres}. But any surjection in $\acat p$ is a composition of small
surjections.
\end{proof}

Let $R$ be any object in $\acat p$. In section \ref{s:ncdefm}, we
described how to, in principle, calculate $\defm(R)$ by
considering the possible left $A$-module structures on the right
$R$-module $(M_i \otimes_k R_{ij})$. The M-free resolutions give
us another way of viewing deformations in $\defm(R)$: By
proposition \ref{p:reslift}, we can view $\defm(R)$ as the set of
liftings of the complex $(F_{\cpd}, d_{\cpd})$ to $R$, up to
equivalence. Using equation \ref{e:diffrep}, each lifting of
complexes corresponds to a family of 1-cochains $\alpha(r_{ij}(l))
\in \hmm^1(L_{\cpd j},L_{\cpd i})$, parametrized by a $k$-basis
for $R$. We leave it as an exercise for the reader to use this
approach to calculate $\defm(R)$ in the case $R = R_{\alpha,
\beta}$ --- this will give a new proof of lemma \ref{l:testdefm}
via the Yoneda representation of $\ext^1_A(M_{\beta},M_{\alpha})$.

\section{Pro-representing hulls of pointed functors}
\label{s:hullf}

We say that a covariant functor $\ff: \acat p \to \sets$ is
\emph{pointed} if $\ff(k^p) = \{ * \}$. In this section, we shall
consider pointed functors defined on the category $\acat p$, and
study their representability. Of course, the motivation for this
is the fact that $\defm$ is such a pointed functor.

Let $R$ be any object of $\cacat p$, and consider the functor
$\h R: \acat p \to \sets$ given by $\h R(S) = \mor(R,S)$ for all
objects $S$ in $\acat p$. The notation $\mor(R,S)$ denotes the
set of morphisms from $R$ to $S$ in the pro-category $\cacat p$.
Then $\h R$ is clearly a pointed functor defined on $\acat p$.

We say that a pointed functor $\ff: \acat p \to \sets$ is
\emph{representable} is $\ff$ is isomorphic to $\h R$ for some
object $R$ in $\acat p$, and pro-representable if $\ff$ is
isomorphic to $\h R$ for some object $R$ in $\cacat p$. However,
it is well-known that deformation functors seldom are
representable or even pro-representable. So a weaker notion is
required, and we shall define the notion of a
\emph{pro-representing hull} of a pointed functor on $\acat p$.
We start by introducing some notation:

Any pointed functor $\ff: \acat p \to \sets$ has an extension
to a functor $\fff : \cacat p \to \sets$ defined on the
pro-category $\cacat p$. This extension is defined by the
formula
    \[ \fff(R) = \lim_{\gets} \; \ff(R/I^n) \]
for any object $R$ in $\cacat p$ with $I=I(R)$. Clearly, any
pointed functor $\ff: \acat p \to \sets$ also has a restriction
to the sub-category $\acat p(n) \subseteq \acat p$ for all $n
\ge 1$. We shall denote this restriction by $\ff_n: \acat p(n)
\to \sets$.

\begin{lem}
\label{l:yoneda}
Let $R$ be an object in $\cacat p$, and let $\ff: \acat p \to
\sets$ be a pointed functor. Then there is a natural
isomorphism of sets $\alpha: \fff (R) \to \mor ( \h R , \ff )$.
\end{lem}
\begin{proof} \noindent
Let $\xi \in \fff(R)$, then $\xi = (\xi_n)$ with $\xi_n \in
\ff(R/I^n)$ for all $n\ge 1$. For any object $S$ in $\acat p$,
we construct a map of sets $\alpha(\xi)_S: \mor(R,S) \to \ff(S)$:
Let $u:R \to S$ be a morphism in $\cacat p$, then $u(I(R))
\subseteq I(S)$, and $I(S)$ is nilpotent since $S$ is in $\acat
r$, so there exists $n\ge 1$ such that $u$ factorizes through
$u_n: R/I(R)^n \to S$. We define $\alpha(\xi)_S(u) =
\ff(u_n)(\xi_n)$, and a straight-forward calculation shows that
this expression is independent upon the choice of $n$, and gives
rise to a natural transformation of functors. Conversely, let
$\phi: \h R \to F$ be a natural transformation of functors
on $\acat p$. Then we define $\xi_n \in \ff(R/I(R)^n)$ to be
$\xi_n = \phi_{R/I(R)^n}(R \to R/I(R)^n)$, where $R \to R/I(R)^n$
is the natural morphism. Again, a straight-forward calculation
shows that $\xi = (\xi_n)$ defines an element in $\fff (R)$, and
that this map of sets defines an inverse to $\alpha$.
\end{proof}

There is also a version of lemma \ref{l:yoneda} for the category
$\acat p(n)$: For an object $R$ in $\acat p(n)$, and a pointed
functor $\ff: \acat p (n) \to \sets$, there is a natural
isomorphism of sets $\alpha_n: \ff(R) \to \mor(\h R, \ff)$. The
construction of this isomorphism is similar to the construction
in lemma \ref{l:yoneda}.

We recall that a morphism $\phi: \ff \to \fg$ of pointed functors
$\ff, \fg: \acat p \to \sets$ is \emph{smooth} if the following
condition holds: For all surjective morphisms $u: R \to S$ in
$\acat p$, the natural map of sets
\begin{equation}
\label{f:smdef}
    \ff(R) \to \ff(S) \bigtimes_{\fg(S)} \fg(R),
\end{equation}
given by $x \mapsto (\ff(u)(x),\phi_R(x))$ for all $x\in \ff(R)$,
is a surjection. Clearly, it is enough to check this for small
surjections in $\acat p$. Also notice that any morphism $\phi:
\ff \to \fg$ of functors naturally extends to a morphism
$\hat{\phi}: \fff \to \ffg$ of functors on $\cacat p$, and if
$\phi$ is a smooth morphism, then $\hat{\phi}_R: \fff(R) \to
\ffg(R)$ is surjective for all objects $R$ in $\cacat p$.

Similarly, we say that a morphism $\phi: \ff \to \fg$ of functors
$\ff,\fg: \acat p(n) \to \sets$ on $\acat p(n)$ is smooth if the
map of sets (\ref{f:smdef}) is surjective for all surjective
morphisms $u:R \to S$ in $\acat p(n)$. Clearly, a morphism $\phi:
\ff \to \fg$ of functors on $\acat p$ is smooth if and only if
the restriction $\phi_n: \ff_n \to \fg_n$ is smooth for all $n
\ge 1$.

Let $\ff$ be a pointed functor on $\acat p$. A \emph{pro-couple}
for $\ff$ is a pair $(R,\xi)$, where $R$ is an object in $\cacat
p$ and $\xi \in \fff(R)$. A morphism $u: (R,\xi) \to (R',\xi')$
of pro-couples is a morphism $u: R \to R'$ in $\cacat p$ such
that $\fff(u)(\xi) = \xi'$. If $(R,\xi)$ is a pro-couple for
$\ff$ such that $R$ is also an object of $\acat p$, then it is
called a \emph{couple} for $\ff$.

We say that a pro-couple $(R,\xi)$ \emph{pro-represents} $\ff$
if $\alpha(\xi): \h R \to \ff$ is an isomorphism of functors on
$\acat p$. If $(R,\xi)$ pro-represents $\ff$ and $(R,\xi)$ is
also a couple for $\ff$, then we say that $(R,\xi)$
\emph{represents} $\ff$. It is clear that if the couple $(R,\xi)$
represents $\ff$, then $(R,\xi)$ is unique up to a unique
isomorphism of couples.

Similarly, let $\ff$ be a pointed functor on $\acat p(n)$. A
couple for $\ff$ is a pair $(R,\xi)$, where $R$ is an object of
$\acat p(n)$ and $\xi \in \ff(R)$. We say that the couple
$(R,\xi)$ represents $\ff$ if and only if $\alpha_n(\xi)$ is an
isomorphism of functors defined on $\acat p(n)$. It is clear
that if this is the case, the couple $(R,\xi)$ is unique up to
a unique isomorphism of couples.

Let $\ff$ be a functor on $\acat p$, and let $(R,\xi)$ be a
pro-couple for $\ff$. For all $n \ge 1$, let $(R_n,\xi_n)$ be
given by $R_n = R/I(R)^n$ and $\xi_n = \ff(u_n)(\xi)$, where
$u_n: R \to R_n$ is the natural surjection. Then $(R_n,\xi_n)$
is a couple for the restriction $\ff_n: \acat p(n) \to \sets$
of $\ff$ for all $n \ge 1$. Notice that $\alpha_n(\xi_n)$ is
the restriction of the morphism $\alpha(\xi)$ to $\acat p(n)$
for all $n \ge 1$. Consequently, $(R,\xi)$ pro-represents $\ff$
if and only if $(R_n,\xi_n)$ represents $\ff_n$ for all $n \ge
1$. In particular, it follows that if $(R,\xi)$ pro-represents
$\ff$, then $(R,\xi)$ is unique up to a unique isomorphism of
pro-couples.

Let $\ff: \acat p \to \sets$ be a pointed functor on $\acat p$.
A \emph{pro-representing hull} of $\ff$ is a pro-couple
$(R,\xi)$ of $\ff$ such that the following conditions hold:
\begin{enumerate}
\item $\alpha(\xi): \h R \to \ff$ is a smooth morphism of functors
on $\acat p$ \item $\alpha_2(\xi_2): \h{R_2} \to \ff_2$ is an
isomorphism of functors on $\acat p(2)$
\end{enumerate}
To simplify notation, we sometimes call the pro-representing
hull $(R,\xi)$ a \emph{hull} of $\ff$.

\begin{prop}
\label{p:uniqueh}
Let $\ff: \acat p \to \sets$ be a pointed functor on $\acat p$,
and assume that $(R,\xi), (R',\xi')$ are pro-representing hulls
of $\ff$. Then there exists an isomorphism of pro-couples $u:
(R,\xi) \to (R',\xi')$.
\end{prop}
\begin{proof} \noindent
Let $\phi=\alpha(\xi), \phi'=\alpha(\xi')$. Since $\phi,\phi'$ are
smooth morphisms, we have that $\phi_{R'}$ and $\phi'_R$ are
surjective. So we can find morphisms $u:(R,\xi) \to (R',\xi')$ and
$v:(R',\xi') \to (R,\xi)$ of pro-couples of $\ff$. The restriction to
$\acat p(2)$ gives us morphisms $u_2: (R_2,\xi_2) \to (R'_2,\xi'_2)$
and $v_2: (R'_2,\xi'_2) \to (R_2,\xi_2)$. But both $(R_2,\xi_2)$ and
$(R'_2,\xi'_2)$ represent $\ff_2$, so $u_2$ and $v_2$ are inverses.
In particular, $\gr_1(u_2)$ and $\gr_1(v_2)$ are inverses, and
$(v\circ u)_2 = v_2 \circ u_2 = \id$. From the proof of proposition
\ref{p:surj}, we see that $\gr (v\circ u)$ is surjective. This means
that $\gr_n(v\circ u)$ is a surjective endomorphims of a finite
dimensional $k$-vector space for all $n\ge 1$, so $\gr (v\circ u)$ is
an isomorphism. By proposition \ref{p:surj}, $v \circ u$ is an
isomorphism as well, and the same holds for $u \circ v$ by a symmetric
argument. It follows that $u$ and $v$ are isomorphisms.
\end{proof}

So if there exists a pro-representing hull of a pointed functor
$\ff$, we know that it is unique, and we shall denote it by
$(H,\xi)$. Notice that $(H,\xi)$ is only unique up to
\emph{non-canonical} isomorphism. By abuse of language, we shall
sometimes omit $\xi$ from the notation, and say that $H$ is the
hull of $\ff$.

\section{Hulls of noncommutative deformation functors}
\label{s:hulldef}

We recall that $k$ is an algebraically closed (commutative) field,
$A$ is an associative $k$-algebra, and $\mfam =\{ M_1, \dots, M_p \}$
is a finite family of left $A$-modules. In this section, we prove
that if the family $\mfam$ satisfy the finiteness condition
(\ref{e:fc}), then there exists a hull $H = H(\mfam)$ of the
noncommutative deformation functor $\defm$. The proof follows Laudal
\cite{lau02}, and the essential point is the following obstruction
calculus:

\begin{prop}
\label{p:obstr}
Let $u: R \to S$ be a small surjective morphism in $\acat p$ with
kernel $K=\ker (u)$, and let $M_S \in \defm(S)$ be a deformation.
Then there exists a canonical obstruction
\[ \obstr \in (\ext^2_A(M_j,M_i) \otimes_k K_{ij} ), \]
such that $\obstr = 0$ if and only if there exists a deformation
$M_R \in \defm(R)$ lifting $M_S$. If this is the case, the set of
deformations in $\defm(R)$ lifting $M_S$ is a torsor under the
$k$-vector space $(\ext_A^1(M_j,M_i) \otimes_k K_{ij} )$.
\end{prop}
\begin{proof} \noindent
We recall from section \ref{s:ncdefm} that up to equivalence,
we may assume that $M_S$ has the following form: $M_S=(M_i
\otimes_k S_{ij})$ with right $S$-module structure given by the
multiplication in $S$, and with left $A$-module structure given
by $k$-linear homomorphisms $a_i: M_i \to \osum (M_j \otimes_k
S_{ji})$ for all $a \in A$. Via the natural projections, the map
$a_i$ gives rise to $k$-linear maps $a_{ji}: M_i \to M_j \otimes_k
S_{ji}$ for $a \in A, \; 1\le i,j \le p$. Since $u$ is surjective,
we may choose $k$-linear maps $L(a)_{ji}: M_i \to M_j \otimes_k
R_{ji}$ such that $(id \otimes u) \circ L(a)_{ji} = a_{ji}$ for
$a \in A, \; 1\le i,j \le p$. Let
\[ L(a) = (L(a)_{ij}) \in (\hmm_k(M_j,M_i \otimes_k R_{ij})), \]
this defines a $k$-linear left action of $A$ on $M_R=(M_i
\otimes_k R_{ij})$, lifting the left $A$-module structure on $M_S$.
We let $Q'= (\hmm_k(M_j,M_i \otimes_k R_{ij}))$, and remark that
this is an associative $k$-algebra in a natural way: We compose
the $k$-linear morphisms in $Q'$ by using the multiplication in
$R$.

For $a,b \in A$, consider the expression $L(ab)-L(a)L(b) \in Q'$.
By the associativity of the left $A$-module structure on $M_S$,
we see that $L(ab)-L(a)L(b) \in Q$, where $Q=(\hmm_k(M_j,M_i
\otimes_k K_{ij})) \subseteq Q'$. Furthermore, we notice that
$Q \subseteq Q'$ is an ideal, and $Q$ has a natural structure as
an $A$-$A$ bimodule via $L$, since $K^2=0$. We define $\psi \in
\hmm_k(A\otimes_k A, Q)$ to be given by $\psi(a,b)= L(ab) -
L(a) L(b)$ for all $a,b \in A$. A straight-forward calculation
shows that $\psi$ is a $2$-cocycle in $\hc^*(A,Q)$, so $\psi$
gives rise to an element $\obstr \in \hh^2(A,Q)$ --- see appendix
\ref{a:hh} for the definition of the Hochschild complex and its
cohomology. Since $K^2=0$, it follows that if $L'$ is another
$k$-linear lifting of the left $A$-module structure on $M_S$, then
the $A$-$A$ bimodule structures of $Q$ given by $L$ and $L'$
coincide. Therefore, $\hh^*(A,Q)$ is independent upon the choice
of $L$, and a straight-forward calculation shows that the same
holds for the obstruction $\obstr$.

We remark that there exists a deformation $M_R \in \defm(R)$
lifting $M_S$ if and only if there exists some $k$-linear lifting
$L': A \to Q'$ of the left $A$-module structure of $M_S$ such
that $L'(ab)=L'(a) L'(b)$ for all $a,b \in A$. Let $\tau = L' -
L$, then $\tau: A \to Q$ is a $k$-linear map, and a straight-forward
calculation shows that $L'(ab) = L'(a) L'(b)$ if and only if the
relation
\[ L(ab) - L(a) L(b) = L(a) \tau(b) - \tau(ab) + \tau(a)
L(b) + \tau(a) \tau(b) \]
holds. Since $K^2=0$, the last term vanishes. The fact that the
above relation holds for all $a,b \in A$ is therefore equivalent
to the fact that $\obstr =0$ in $\hh^2(A,Q)$. So we have established
that there exists a canonical obstruction $\obstr \in \hh^2(A,Q)$
such that $\obstr=0$ if and only if there is a lifting of $M_S$ to
$R$.

Assume that $L: A \to Q'$ is such that $L(ab) = L(a) L(b)$ for all
$a,b \in A$, that is, such that it defines a deformation $M_R$ lying
over $M_S$. For any other $k$-linear lifting $L': A \to Q'$ of the left
$A$-module structure on $M_S$, we may consider the difference $\tau =
L' - L : A \to Q$. A straight-forward calculation shows that $\tau$ is
a $1$-cocycle in $\hc^*(A,Q)$ if and only if $L'(ab)=L'(a) L'(b)$ for
all $a,b \in A$, that is, if and only if $L'$ defines a left $A$-module
structure on $M_R$. Furthermore, we have that $L$ and $L'$ give rise to
equivalent deformations if and only if $\tau$ is a $1$-coboundary: It is
clear that any equivalence between the left $A$-module structures of
$M_R=(M_i \otimes_k R_{ij})$ given by $L$ and $L'$ has the form $id+\psi$,
where $\psi\in Q$. Furthermore, the map $id+\psi: M_R \to M_R$ (with the
left $A$-module structure from $L'$ and $L$ respectively) is a left
$A$-module homomorphism if and only if $L(r)(id + \psi) = (id + \psi) L'(r)$
holds for all $a \in A$, and this last condition is equivalent with the fact
that $\tau=d(\psi)$, so that $\tau$ is a $1$-coboundary. If $\tau$ is a
$1$-boundary in $\hc^*(A,Q)$, it is also clear that $id+\psi$ defines an
equivalence between the two deformations given by $L$ and $L'$. Therefore,
the set of deformations $M_R$ lying over $M_S$ is a torsor under the
$k$-vector space $\hh^1(A,Q)$.

To end the proof, we have to show that there are isomorphisms of
$k$-vector spaces $\hh^n(A,Q) \cong (\ext^n_A(M_j,M_i) \otimes_k
K_{ij})$ for $n=1,2$: Since $L(a)$ is a lifting to $M_R$ of the
left multiplication of $a$ on $M_S$ (satisfying equation \ref{f:fdef}),
$L(a)$ satisfies equation \ref{f:fdef} as well. That is, we have
$L(a)_{ji}(m_i)-\delta_{ij}(am_i)\otimes e_i \in M_j \otimes_k I_{ji}$
for all $a \in A, \; m_i \in M_i, \; 1 \le i,j \le p$. Since $K^2=0$, this
means that the $A$-$A$ bimodule structure of $Q$ defined via $L$
coincides with the following natural one: Since $M_i, M_j \otimes_k
K_{ji}$ are left $A$-modules, we have that $Q_{ij}=\hmm_k(M_j,M_i
\otimes_k K_{ij})$ and $Q = \oplus Q_{ij}$ has natural $A$-$A$ bimodule
structures. Clearly, we have
\[ \hh^n(A,Q) \cong \osum_{i,j} \hh^n(A,Q_{ij}) = (\hh^n(A,Q_{ij})). \]
By appendix \ref{a:hh}, proposition \ref{p:hhext}, we have that
$\hh^n(A,Q_{ij}) \cong \ext^n_A(M_j,M_i \otimes_k K_{ij})$ for $n \ge 0$.
Moreover, $\ext^n_A(M_j,M_i \otimes_k K_{ij}) \cong \ext^n_A(M_j,M_i)
\otimes_k K_{ij}$ since $K_{ij}$ is a $k$-vector space of finite
dimension. This completes the proof of the proposition.
\end{proof}

We remark that it is easy to find an alternative proof of proposition
\ref{p:obstr} using resolutions and the Yoneda representation of
$\ext^n_A(M_i,M_j)$. This is straight-forward, but makes essential use
of proposition \ref{p:reslift}.

Also notice that the obstruction calculus is functorial in the following
sense: Let $u: R \to S$ and $u':R' \to S'$ be two small surjections in
$\acat p$, and write $K=\ker(u)$ and $K'=\ker(u')$. Assume that $v:R \to
R'$ and $w:S \to S'$ are morphisms such that $u' \circ v = w \circ u$.
Then $v(K) \subseteq K'$, and the map $v$ induces a $k$-linear map of
obstruction spaces
\[ (\ext^2_A(M_j,M_i) \otimes_k K_{ij}) \to (\ext^2_A(M_j,M_i) \otimes_k
K'_{ij}). \]
If $M_S$ is a deformation of $\mfam$ to $S$ and $M_{S'}=\defm(w)(M_S)$ is
the corresponding deformation to $S'$, then this map of obstruction
spaces maps $o(u,M_S)$ to $o(u',M_{S'})$. This follows from the proof of
proposition \ref{p:obstr}.

Let us start the construction of the pro-representing hull $(H,\xi)$
of $\defm$, using the obstruction calculus for $\defm$ given above. From
now on, we shall assume that the family $\mfam$ satisfy the finiteness
condition
\begin{equation}
\tag{\bfseries FC}
\dim_k \ext^n_A(M_i,M_j) \text{ is finite for } 1 \le i,j \le p, \;
n=1,2.
\end{equation}
We fix the following notation: Let $\{x_{ij}(l): 1 \le l \le d_{ij}
\}$ be a basis for $\ext^1_A(M_j,M_i)^*$ and let $\{y_{ij}(l): 1 \le l
\le r_{ij} \}$ be a basis for $\ext^2_A(M_j,M_i)^*$ for $1 \le i,j \le
p$, with $d_{ij} = \dim_k \ext^1_A(M_j,M_i)$ and $r_{ij} = \dim_k
\ext^2_A(M_j,M_i)$. Moreover, we consider the formal matrix rings in
$\cacat p$ corresponding to these vector spaces, and denote them by
$\infdef = \hat{R}(\{\ext^1_A(M_j,M_i)^*\})$ and $\obstrdef =
\hat{R}(\{\ext^2_A(M_j,M_i)^*\})$.

First, let us show that $\defm$ restricted to $\acat p(2)$ is representable:
We define $H_2$ to be the object $H_2=\infdef_2=\infdef / I(\infdef)^2$ in
$\acat p(2)$. For all objects $R$ in $\acat p(2)$, we get $\mor (H_2,R)
\cong (\hmm_k(\ext^1_A(M_j,M_i)^*,I(R)_{ij})) \cong (\ext^1_A(M_j,M_i)
\otimes_k I(R)_{ij})$, and $\defm(R) \cong (\ext^1_A(M_j,M_i) \otimes_k
I(R)_{ij})$ by proposition \ref{p:obstr} applied to the small surjection $R
\to k^p$. The isomorphisms we obtain in this way are compatible, so they
induce an isomorphism $\phi_2: \h {H_2} \to \defm$ of functors on $\acat
p(2)$. From the version of lemma \ref{l:yoneda} for the category $\acat
p(2)$, we see that there is a unique deformation $\xi_2 \in \defm(H_2)$
such that $\alpha_2(\xi_2)=\phi_2$. By definition, $(H_2,\xi_2)$ represents
the deformation functor $\defm$ restricted to $\acat p(2)$.

Let us also give an explicit description of the deformation $\xi_2$: We
have $H_2=\infdef_2$, so let us denote by $\epsilon_{ij}(l)$ the image of
$x_{ij}(l)$ in $H_2$ for $1 \le i,j \le p, \; 1 \le l \le d_{ij}$. In this
notation, $\xi_2$ is represented by the right $H_2$-module $(M_i \otimes_k
(H_2)_{ij})$, with left $A$-module structure defined by
\[ a (m_j \otimes e_j) = am_j \otimes e_j + \sum_{i,l} \psi_{ij}^l(a)(m_j)
\otimes \epsilon_{ij}(l) \]
for all $a \in A, \; m_j \in M_j, \; 1 \le j \le p$, where $\psi_{ij}^l
\in \der_k(A,\hmm_k(M_j,M_i))$ is a representative of $x_{ij}(l)^* \in
\ext^1_A(M_j,M_i)$ via Hochschild cohomology.

There is also an alternative description of $\xi_2$ using M-free
resolutions and the Yoneda representation of $\ext^1_A(M_i,M_j)$:
Let $\alpha(\epsilon_{ij}(l)) \in \hmm^1(L_{\cpd j},L_{\cpd i})$
be a $1$-cocycle representing $x_{ij}(l)^* \in \ext^1_A(M_j,M_i)$
for $1 \le i,j \le p, \; 1 \le l \le d_{ij}$. Then by
construction, the formula
    \[ d^{H_2}_m = \sum_i d_{m,i} \otimes e_i + \sum_{i,j,l}
    \alpha(\epsilon_{ij}(l))_m \otimes \epsilon_{ij}(l) \]
defines a differential which lifts the complex
$(F_{\cpd},d_{\cpd})$ to $H_2$. By proposition \ref{p:reslift},
the lifted complex is in fact an M-free resolution of some
deformation of $\mfam$ to $H_2$, and this deformation is $\xi_2
\in
\defm(H_2)$.

\begin{thm}
\label{t:exh}
Assume that $\dim_k \ext^n_A(M_i,M_j)$ is finite for $1 \le i,j \le p, \;
n=1,2$. Then there exists a morphism $o:\obstrdef \to \infdef$ in $\cacat
p$ such that $H(\mfam)= \infdef \hat{\otimes}_{\obstrdef} k^p$ is a
pro-representing hull for $\defm$.
\end{thm}
\begin{proof} \noindent
For simplicity, let us write $I$ for the ideal $I=I(\infdef)$, and for all
$n \ge 1$, let us write $\infdef_n$ for the quotient $\infdef_n=\infdef/I^n$,
and $t_n: \infdef_{n+1} \to \infdef_n$ for the natural morphism. From the
paragraphs preceding this theorem, we know that $(H_2,\xi_2)$ represents
$\defm$ restricted to $\acat p(2)$. Let $o_2: \obstrdef \to \infdef_2$ be
the trivial morphism given by $o_2(I(\obstrdef))=0$ and let $a_2 = I^2$,
then $H_2 = \infdef / a_2 \cong \infdef_2 \otimes_{\obstrdef} k^p$. Using
$o_2$ and $\xi_2$ as a starting point, we shall construct $o_{n}$ and $\xi_n$
for $n \ge 3$ by an inductive process. So let $n \ge 2$, and assume that the
morphism $o_n: \obstrdef \to \infdef_n$ and the deformation $\xi_n \in
\defm(H_n)$ is given, with $H_n = \infdef_n \otimes_{\obstrdef} k^p$. We shall
also assume that $t_{n-1} \circ o_n = o_{n-1}$ and that $\xi_n$ is a lifting
of $\xi_{n-1}$.

Let us now construct the morphism $o_{n+1}: \obstrdef \to
\infdef_{n+1}$: We let $a'_n$ be the ideal in $\infdef_n$
generated by $o_n(I(\obstrdef))$. Then $a'_n = a_n/I^n$ for an
ideal $a_n \subseteq \infdef$ with $I^n \subseteq a_n$, and $H_n
\cong \infdef / a_n$. Let $b_n = I a_n + a_n I$, then we obtain
the following commutative diagram:
\[
\xymatrix{
{\obstrdef} \ar[rd]^{o_n} & {\infdef_{n+1}} \ar[d] \ar[r] & {\infdef/b_n}
\ar[d] \\
& {\infdef_n} \ar[r] & H_n=\infdef/a_n,
}
\]
Observe that $\infdef/b_n \to \infdef/a_n$ is a small surjection.
So by proposition \ref{p:obstr}, there is an obstruction
$o'_{n+1}=o(\infdef/b_n \to H_n,\xi_n)$ for lifting $\xi_n$ to
$\infdef/b_n$, and we have
\[ o'_{n+1} \in (\ext^2_A(M_j,M_i)\otimes_k (a_n/b_n)_{ij}) \cong
(\hmm_k( \gr_1(\obstrdef)_{ij} , (a_n/b_n)_{ij})). \]
Consequently, we obtain a morphism $o'_{n+1}: \obstrdef \to \infdef/b_n$. Let
$a''_{n+1}$ be the ideal in $\infdef/b_n$ generated by
$o'_{n+1}(I(\obstrdef))$. Then $a''_{n+1}=a_{n+1}/b_n$ for an ideal $a_{n+1}
\subseteq \infdef$ with $b_n \subseteq a_{n+1} \subseteq a_n$. We define
$H_{n+1}=\infdef/a_{n+1}$ and obtain the following commutative diagram:
\[
\xymatrix{
{\obstrdef} \ar[rd]_{o_n} \ar@/^1pc/[rr]^{o'_{n+1}}
& {\infdef_{n+1}} \ar[d] \ar[r] & {\infdef/b_n} \ar[r] \ar[d] & H_{n+1}=
\infdef/a_{n+1} \ar[dl] \\
& {\infdef_n} \ar[r] &
H_n=\infdef/a_n & }
\]
By the choice of $a_{n+1}$, the obstruction for lifting $\xi_n$ to
$H_{n+1}$ is zero. We can therefore find a lifting $\xi_{n+1} \in
\defm(H_{n+1})$ of $\xi_n$ to $H_{n+1}$.

The next step of the construction is to find a morphism $o_{n+1}:\obstrdef
\to \infdef_{n+1}$ which commutes with $o'_{n+1}$ and $o_n$: We know that
$t_{n-1} \circ o_n = o_{n-1}$, which means that $a_{n-1} = I^{n-1} + a_n$.
For simplicity, let us write $O(K)=(\hmm_k(\gr_1(\obstrdef)_{ij},K_{ij}))$
for any ideal $K \subseteq \infdef$. Consider the following commutative
diagram of $k$-vector spaces, in which the columns are exact:
\[
\xymatrix{
{0} \ar[d] & {0} \ar[d] \\
{O(b_n/I^{n+1})} \ar[d] \ar[r]^{j_n} & {O(b_{n-1}/I^n)} \ar[d] \\
{O(a_n/I^{n+1})} \ar[d]^{r_{n+1}} \ar[r]^{k_n} & {O(a_{n-1}/I^n)}
\ar[d]^{r_n} \\
{O(a_n/b_n)} \ar[d] \ar[r]^{l_n} & {O(a_{n-1}/b_{n-1})} \ar[d] \\
{0} & {0} }
\]
We may consider consider $o_n$ as an element in $O(a_{n-1}/I^n)$, since
$a_n \subseteq a_{n-1}$. On the other hand, $o'_{n+1} \in O(a_n/b_n)$.
Let $o'_n = r_n(o_n)$, then the natural map $\infdef/b_n \to \infdef/
b_{n-1}$ maps the obstruction $o'_{n+1}$ to the obstruction $o'_n$ by
the second remark following proposition \ref{p:obstr}. This implies that
$o'_{n+1}$ commutes with $o'_n$, so $l_n(o'_{n+1}) = o'_n = r_n(o_n)$.
But we have $o_n(I(\obstrdef)) \subseteq a_n$, so we can find an element
$\overline{o}_{n+1} \in O(a_n/I^{n+1})$ such that
$k_n(\overline{o}_{n+1})=o_n$. Since $a_{n-1}=a_n + I^{n-1}$, $j_n$ is
surjective. Elementary diagram chasing using the snake lemma implies
that we can find $o_{n+1} \in O(a_n/I^{n+1})$ such that $r_{n+1}(o_{n+1})
= o'_{n+1}$ and $k_n(o_{n+1}) = o_n$. It follows that the obstruction
$o_{n+1}$ defines a morphism $o_{n+1}: \obstrdef \to \infdef_{n+1}$
compatible with $o_n$ such that $\infdef_{n+1} \otimes_{\obstrdef} k^p
\cong H_{n+1}$.

By induction, it follows that we can find a morphism $o_n: \obstrdef \to
\infdef_n$ and a deformation $\xi_n \in \defm(H_n)$, with $H_n=\infdef_n
\otimes_{\obstrdef} k^p$, for all $n\ge 1$. From the construction, we see
that $t_{n-1} \circ o_n = o_{n-1}$ for all $n \ge 2$, so we obtain a
morphism $o: \obstrdef \to \infdef$ by the universal property of the
projective limit. Moreover, the induced morphisms $h_n: H_{n+1} \to H_n$
are such that  $\xi_{n+1} \in \defm(H_{n+1})$ is a lifting of $\xi_n \in
\defm(H_n)$ to $H_{n+1}$. Notice that $I(H_n)^n = 0$ and that
$H_n/I(H_n)^{n-1} \cong H_{n-1}$ for all $n \ge 2$. It follows that
$H/I(H)^n = H_n$ for all $n \ge 1$, so $H$ is an object of the
pro-category $\cacat p$. Let $\xi = (\xi_n)$, then clearly $\xi \in
\defm(H)$, so $(H,\xi)$ is a pro-couple for $\defm$. It remains to show
that $(H,\xi)$ is a pro-representable hull for $\defm$.

It is clearly enough to show that $(H_n,\xi_n)$ is a pro-representing hull
for $\defm$ restricted to $\acat p(n)$ for all $n\ge 3$. So let $\phi_n =
\alpha_n(\xi_n)$ be the morphism of functors on $\acat p(n)$ corresponding
to $\xi_n$. We shall prove that $\phi_n$ is a smooth morphism. So let $u:R
\to S$ be a small surjection in $\acat p(n)$, and assume that $M_R \in
\defm(R)$ and $v\in \mor(H_n,S)$ are given such that $\defm(u)(M_R)=
\defm(v)(\xi_n)=M_S$. Let us consider the following commutative diagram:
\[
\xymatrix{
& {\infdef/b_n} \ar[d] & \\
{\infdef} \ar[ru] \ar[r] \ar[rd] & {H_{n+1}} \ar[d] &
{R} \ar[d]^u \\
& {H_n} \ar[r]_v & {S}
}
\]
Let $v': \infdef \to R$ be any morphism making the diagram commutative.
Then $v'(a_n) \subseteq K$, where $K=\ker(u)$, so $v'(b_n)=0$. But the
induced map $\infdef/b_n \to R$ maps the obstruction $o'_{n+1}$ to
$o(u,M_S)$, and we know that $o(u,M_S)=0$. So we have $v'(a_{n+1})=0$,
and $v'$ induces a morphism $v': H_{n+1} \to R$ making the diagram
commutative. Since $v'(I(H_{n+1})^n)=0$, we may consider $v'$ a map from
$H_{n+1}/I(H_{n+1})^n \cong H_n$. So we have constructed a map $v' \in
\mor(H_n,R)$ such that $u \circ v' = v$. Let $M'_R = \defm(v')(\xi_n)$,
then $M'_R$ is a lifting of $M_S$ to $R$. By proposition \ref{p:obstr},
the difference between $M_R$ and $M'_R$ is given by an element
\[ d \in (\ext^1_A(M_j,M_i) \otimes_k K_{ij}) =
(\hmm_k(\gr_1(\infdef)_{ij},K_{ij})). \] Let $v'': \infdef \to R$
be the morphism given by $v''(x_{ij}(l)) = v'(x_{ij}(l)) +
d(\overline{x_{ij}(l)})$ for $1 \le i,j \le p, \; 1 \le l \le
d_{ij}$. Since $a_{n+1} \subseteq I(\infdef)^2$, we have
\[ v''(a_{n+1}) \subseteq v'(a_{n+1}) + I(R)K + KI(R) + K^2. \]
But $u$ is small, so $v''(a_{n+1})=0$ and $v''$ induces a morphism $v'':
H_n \to R$. Clearly, $u \circ v'' = u \circ v' = v$, and $\defm(v'')
(\xi_n) = M_R$ by construction. It follows that $\phi_n$ is smooth for
all $n \ge 3$.
\end{proof}

We remark that the conclusion of the theorem still holds if we relax the
finiteness condition (\ref{e:fc}). If we only assume that
\[ \dim_k \ext^1_A(M_i,M_j) \text{ is finite for } 1 \le i,j \le p, \]
then the object $\obstrdef$ is in $\Acat p$, but not necessarily in
$\cacat p$. However, the rest of the proof is still valid as stated, so
the finiteness condition on $\ext^2_A(M_i,M_j)$ is clearly not essential.

In general, it is possible to generalize theorem \ref{t:exh} to the case
when $\ext^n_A(M_i,M_j)$ has countable dimension as a vector space over
$k$ for $1 \le i,j \le p, \; n=1,2$, see Laudal \cite{lau79}. However,
we shall always assume (\ref{e:fc}) in this paper.

Assume that $\mfam$ satisfy (\ref{e:fc}). If $\ext^2_A(M_i,M_j)=0$ for
$1 \le i,j \le p$, we say that the deformation functor $\defm$ is
\emph{unobstructed}. For instance, $\defm$ is unobstructed for any
finite family $\mfam$ of left $A$-modules satisfying (\ref{e:fc}) if
$A$ is left hereditary (that is, the left global homological dimension
of $A$ is at most $1$). If $\defm$ is unobstructed, $H = \infdef$ is
the hull of $\defm$.

In general, $\defm$ can be obstructed, and there is no simple formula
for the hull $H$ of $\defm$ if this is the case. However, there exists
an algorithm for calculating the hull $H$ using matric Massey products.
In the next sections, we shall introduce the matric Massey products
and explain how the hull can be calculated when $\mfam$ satisfy
(\ref{e:fc}).

\section{Immediately defined matric Massey products}
\label{s:mmp}

We recall that $k$ is an algebraically closed (commutative) field,
$A$ is an associative $k$-algebra, and $\mfam =\{ M_1, \dots, M_p
\}$ is a finite family of left $A$-modules. From now on, we also
assume that the family $\mfam$ satisfy the finiteness condition
(\ref{e:fc}). In this section, we shall define the immediately
defined matric Massey products and their defining system, and show
how to calculate these products using matrices.

Let us fix a monomial $X \in I(\infdef)$ of type $(i,j)$ and degree
$n \ge 2$. Then we can write $X$ uniquely in the form
    \[ X = x_{i_0 i_1}(l_1) \; x_{i_1 i_2}(l_2) \;  \dots \;
    x_{i_{n-1} i_n}(l_n), \]
where $(i_0,i_n) = (i,j)$. Let $X'$ be another monomial in $\infdef$.
We shall say that $X'$ \emph{divides} $X$ if there exist monomials
$X(l), X(r) \in \infdef$ such that $X = X(l) X' X(r)$, and write
$X' \dvd X$ if this is the case.

Consider the set of monomials $\{ X' \in I(\infdef): X' \ndv X \}$,
and denote by $J(X)$ the ideal in $\infdef$ generated by these
monomials. We define $R(X) = \infdef/J(X)$ and $S(X) = R(X) / (X) =
\infdef/(J(X),X)$. Then the natural map
    \[ \pi(X): R(X) \to S(X) \]
is a small surjection in $\acat p$, and it has a $1$-dimensional
kernel which is generated by the monomial $X$. We write $I(X) =
I(S(X))$ and $S(X)_n = S(X)/I(X)^n$ for all $n \ge 1$.

Let us consider the set $B(X) = \{ (i,j,l): 1 \le i,j \le p, \; 1 \le
l \le d_{ij}, \; x_{ij}(l) \dvd X \}$, and denote by $v_{ij}(l)$ the
image of $x_{ij}(l)$ in $S(X)_2$ for all $(i,j,l) \in B(X)$. Then the
set
    \[ \{ v_{ij}(l) : (i,j,l) \in B(X) \} \]
is a natural $k$-basis for $I(X)/I(X)^2$.

Assume that a morphism $\phi(X): H \to S(X)$ is given, and denote the
composition of $\phi(X)$ with the natural morphism $S(X) \to S(X)_2$
by $\phi(X)_2: H \to S(X)_2$. This morphism can be written uniquely
in the form
    \[ \phi(X)_2 = \sum_{(i,j,l) \in B(X)} \alpha_{ij}(l) \otimes
    v_{ij}(l), \]
where $\alpha_{ij}(l) \in \ext^1_A(M_j,M_i)$ for all $(i,j,l) \in
B(X)$.

Conversely, consider a family $\{ \alpha_{ij}(l) \in \ext^1_A(M_j,
M_i): (i,j,l) \in B(X) \}$ of extensions indexed by $B(X)$,
corresponding to a morphism $\phi(X)_2: H \to S(X)_2$ given by
$\phi(X)_2 = \sum \alpha_{ij}(l) \otimes v_{ij}(l)$. If there exists
a lifting of $\phi(X)_2$ to a morphism $\phi(X): H \to S(X)$, we say
that the \emph{matric Massey product}
        \[ \langle \ul \alpha ; X \rangle = \langle \alpha_{i_0,i_1}
    (l_1), \alpha_{i_1,i_2}(l_2), \dots, \alpha_{i_{n-1},i_n}
    (l_n) \rangle \]
is defined, and that $\phi(X)$ is a \emph{defining system} for this
matric Massey product. If this is the case, we denote the deformation
induced by the defining system $\phi(X)$ by $M_X \in \defm(S(X))$,
and by proposition \ref{p:obstr}, the obstruction for lifting $M_X$
to $R(X)$ is an element
    \[ o(\pi(X), M_X) \in (\ext^2_A(M_j,M_i) \otimes_k K(X)_{ij})
    \cong \ext^2_A(M_j,M_i), \]
where $K(X) = \ker(\pi(X)) \cong k X$. In general, this element
depends upon the deformation $M_X$, and therefore on the defining
system $\phi(X)$. We define the value of the matric Massey product
to be
    \[ \langle \ul \alpha ; X \rangle = \langle \alpha_{i_0,i_1}
    (l_1), \alpha_{i_1,i_2}(l_2), \dots, \alpha_{i_{n-1},i_n}
    (l_n) \rangle = o(\pi(X), M_X). \]
Consequently, the value of the matric Massey product $\langle \ul
\alpha ; X \rangle$ will in general depend upon the chosen defining
system.

Let us fix the monomial $X$. Then the matric Massey product $\ul
\alpha \mapsto \langle \ul \alpha ; X \rangle$ is a not everywhere
defined $k$-linear map
\[
\xymatrix{
\ext^1_A(M_{i_1},M_{i_0}) \otimes_k \dots \otimes_k
\ext^1_A(M_{i_n},M_{i_{n-1}}) \ar@{-->}[r] & \ext^2_A(M_{i_n},
M_{i_0}). }
\]
In fact, this map is defined for $\ul \alpha$ if and only if the
morphism $\phi(X)_2: H \to S(X)_2$ corresponding to $\ul \alpha$
can be lifted to a morphism $\phi(X): H \to S(X)$. Moreover, even
when this map is defined for $\ul \alpha$, it is not necessarily
uniquely defined: In general, its value $\langle \ul \alpha ; X
\rangle$ depends upon the chosen lifting $\phi(X)$, the defining
system. The matric Massey products $\langle \ul \alpha ; X
\rangle$ defined above are called the \emph{immediately defined
matric Massey products}.

We remark that if $X$ is a monomial of degree $n=2$, then the
situation is much simpler: We have $S(X) = S(X)_2$, so the
matric Massey product $\langle \ul \alpha ; X \rangle$ is
uniquely defined for any family of extensions $\{ \alpha_{ij}(l):
(i,j,l) \in B(X) \}$. In fact, the matric Massey product is just
the usual \emph{cup product} in this case.

Let us fix a monomial $X \in I(\infdef)$ of degree $n \ge 2$.
Then there exists a natural family of extensions indexed by
$B(X)$ given by $\alpha_{ij}(l) = x_{ij}(l)^*$,
    \[ \{ x_{ij}(l)^* \in \ext^1_A(M_j,M_i): (i,j,l) \in
    B(X) \}. \]
The matric Massey products of these extensions are the ones that
we shall use for the construction of the hull $H$ of $\defm$ in
the next section. We therefore introduce the notation
    \[ \langle \ul x^* ; X \rangle = \langle x_{i_0 i_1}
    (l_1)^*, x_{i_1 i_2}(l_2)^*, \dots, x_{i_{n-1} i_n}
    (l_n)^* \rangle \]
for their immediately defined matric Massey products.

The matric Massey products are called \emph{matric} because these
products (and their defining systems) can be described completely
in terms of linear algebra and matrices. We shall end this section
by giving such a description.

Let $\{ \alpha_{ij}(l) \in \ext^1_A(M_j,M_i): (i,j,l) \in B(X) \}$
be a family of extensions indexed by $B(X)$, and consider the
corresponding matric Massey product
\begin{equation} \label{e:mmp}
\langle \ul \alpha ; X \rangle = \langle \alpha_{i_0 i_1}(l_1),
\alpha_{i_1 i_2}(l_2), \dots, \alpha_{i_{n-1} i_n}(l_n) \rangle.
\end{equation}
We assume that there exists a defining system $\phi(X): H \to
S(X)$ for this matric Massey product. Then $\phi(X)$ induces a
deformation $M_X \in \defm(S(X))$. We notice that the matric
Massey product (\ref{e:mmp}) only depends upon this deformation.
By abuse of language, we shall therefore let the notion
\emph{defining system} refer to the deformation $M_X$ as well as
the morphism $\phi(X): H \to S(X)$ which induces $M_X$.

We know that any deformation $M_X \in \defm(S(X))$ can be described
by a complex which lifts $(F_{\cpd},d_{\cpd})$ to $S(X)$. Such a
complex is given by differentials of the form
    \[ d^{S(X)}_m: (L_{m+1,i} \otimes_k S(X)_{ij}) \to
    (L_{m,i} \otimes_k S(X)_{ij}). \]
We write $v(X')$ for the image of $X'$ in $S(X)$ whenever $X'$
is a monomial in $\infdef$, and define $\ol B(X) = \{ X' \in
I(\infdef): X' \text{ is a monomial such that } X' \dvd X \} \cup
\{ e_1, \dots, e_p \}$. Then the set
    \[ \{ v(X'): X' \in \ol B(X) \} \]
is a natural $k$-basis for $S(X)$, and $\ol B(X)$ contains
$\{ x_{ij}(l): (i,j,l) \in B(X) \}$ and $\{ e_1, \dots, e_p \}$ as
subsets. Let us write
$\ol B(X)_{ij} = \ol B(X) \cap S(X)_{ij}$ for $1 \le i,j \le p$.
With this notation, the above differentials have the form
    \[ d^{S(X)}_m = \sum_{1 \le i \le p} d_{m,i} \otimes e_i \; +
    \sum_{X' \in \overline B(X)} \alpha(X')_m \otimes v(X'), \]
where $\alpha(X') \in \hmm^1_A(L_{\cpd j},L_{\cpd i})$ is a
$1$-cochain whenever $X' \in \ol B(X')_{ij}$.

Let $d^{S(X)}_m$ be arbitrary maps between M-free modules over
$S(X)$ defined by a family of $1$-cochains $\{ \alpha(X'): X' \in
\ol B(X) \}$ as above. These maps lifts the complex $(F_{\cpd},
d_{\cpd})$ if $\alpha(e_i)=d_{\cpd i}$ for $1 \le i \le p$.
Moreover, these maps are differentials if and only if the
following condition holds: For all monomials $Z \in \ol B(X)$ and
for all integers $m \ge 0$, we have
\begin{equation}
\sum_{\substack{X', X'' \in \ol B(X) \\ X' X'' = Z }} \alpha(X')_m
\circ \alpha(X'')_{m+1} = \sum_{\substack{X', X'' \in \ol B(X) \\
X' X'' = Z}} \alpha(X'')_{m+1} \alpha(X')_m = 0.
\end{equation}
In the first sum, the symbol $\circ$ denotes composition of maps.
We recall that each of the maps involved can be considered as right
multiplication by a matrix. In the second summation, we identify
the maps with such matrices, and re-write the composition of maps
as multiplication of the corresponding matrices.

Assume that these conditions hold. Then the family $\{ \alpha(X'):
X' \in \ol B(X) \}$ of $1$-cochains defines a lifting of complexes
of $(L_{\cpd},d_{\cpd})$ to $S(X)$ given by the differentials
$d^{S(X)}$ as above, and this lifting corresponds to a deformation
$M_X \in \defm(S(X))$. The deformation $M_X$ is a defining system
for the matric Massey product (\ref{e:mmp}) if and only if
$\alpha(X')$ is a $1$-cocycle which represents $\alpha_{ij}(l) \in
\ext^1_A(M_j,M_i)$ whenever $X' = x_{ij}(l)$ for some $(i,j,l) \in
B(X)$. In this case, we shall refer to the family of $1$-cochains
$\{ \alpha(X'): X' \in \ol B(X) \}$ as a \emph{defining system}
for the matric Massey product (\ref{e:mmp}).

Finally, assume that the family of $1$-cochains $\{ \alpha(X'): X'
\in \ol B(X) \}$ is a defining system of the matric Massey product
(\ref{e:mmp}). Then the value of this matric Massey product is
given by
\begin{equation}
\langle \ul \alpha ; X \rangle_m = \sum_{\substack{X',X'' \in \ol
B(X) \\ X' X'' = X}} \alpha(X'')_{m+1} \alpha(X')_m
\end{equation}
for all $m \ge 0$, where the multiplication denotes matrix
multiplication of the corresponding matrices.

\begin{prop} \label{p:mmp}
Let $\{ \alpha_{ij}(l) \in \ext^1_A(M_j,M_i): (i,j,l) \in B(X) \}$
be a family of extensions. A defining system for the matric Massey
product
    \[ \langle \ul \alpha ; X \rangle = \langle \alpha_{i_0
    i_1}(l_1), \dots, \alpha_{i_{n-1} i_n}(l_n) \rangle \]
corresponds to a family $\{ \alpha(X') \in \hmm^1_A(L_{\cpd j},
L_{\cpd i}): 1 \le i,j \le p, \; X' \in \ol B(X)_{ij} \}$ of
$1$-cochains satisfying the following conditions:
\begin{itemize}
\item $\alpha(e_i)=d_{\cpd i}$ for $1 \le i \le p$,
\item $\alpha(X')$ is a $1$-cocycle representing $\alpha_{ij}(l)$
whenever $X' = x_{ij}(l)$ for some $(i,j,l) \in B(X)$,
\item For all $Z \in \ol B(X)$ and for all $m \ge 0$, we have
\begin{equation*}
\sum_{\substack{X', X'' \in \ol B(X) \\ X' X'' = Z}}
\alpha(X'')_{m+1} \; \alpha(X')_m = 0.
\end{equation*}
\end{itemize}
Moreover, given such a family of $1$-cochains, the matric Massey
product $\langle \ul \alpha ; X \rangle$ is represented by the
$2$-cocyle given by
    \[ \langle \ul \alpha ; X \rangle_m = \sum_{\substack{X',
    X'' \in \ol B(X) \\ X' X'' = X}} \alpha(X'')_{m+1} \;
    \alpha(X')_m \]
for all $m \ge 0$.
\end{prop}

Hence we have described the immediately defined matric Massey
products and their defining systems in terms of linear algebra and
matrices, as we set out to do. We remark that the description
given in proposition \ref{p:mmp} is extremely useful for doing
concrete calculations with matric Massey products, and even for
implementing such computations on computers. It also justifies the
name \emph{matric}.

\section{Calculating hulls using matric Massey products}
\label{s:ch}

We recall that $k$ is an algebraically closed (commutative) field,
$A$ is an associative $k$-algebra, and $\mfam =\{ M_1, \dots, M_p
\}$ is a finite family of left $A$-modules. We also assume that
the family $\mfam$ satisfy the finiteness condition (\ref{e:fc}).
In this section, we show how to calculate the hull $H$ of the
deformation functor $\defm$ using matric Massey products.

By theorem \ref{t:exh}, there exists an obstruction morphism $o:
\obstrdef \to \infdef$ in $\cacat p$ such that $H = \infdef
\hat{\otimes}_{\obstrdef} k^p$ is a hull for the deformation functor
$\defm$. We shall write $I = I(\infdef)$ and $f_{ij}(l) =
o(y_{ij}(l))$ for $1 \le i,j  \le p, \; 1 \le l \le r_{ij}$. Then
$f_{ij}(l)$ is a formal power series in $I^2_{ij}$ by construction.
Let us define $a \subseteq \infdef$ to be the ideal generated by
$\{ f_{ij}(l): 1 \le i,j \le p, \; 1 \le l \le r_{ij} \}$. Then $a
\subseteq I^2$, and we have
\begin{equation*}
H = \infdef \hat{\otimes}_{\obstrdef} k^p \cong \infdef / a.
\end{equation*}
We shall use the matric Massey products from section \ref{s:mmp} to
calculate the coefficients of the power series $f_{ij}(l)$. Clearly,
this is sufficient to determine the hull $H$.

Let us fix an integer $N \ge 2$ such that $a \subseteq I^N$. This
is always possible, since $a \subseteq I^2$. So $f_{ij}(l) \in I^N$
for all $f_{ij}(l)$, and we can write  $f_{ij}(l)$ in the form
    \[ f_{ij}(l) = \sum_{\ord X = N} a_{ij}^l(X) \cdot X +
    \sum_{\ord X > N} a_{ij}^l(X) \cdot X \]
for $1 \le i,j \le p, \; 1 \le l \le r_{ij}$, with $a_{ij}^l(X) \in
k$ for all monomials $X \in I^N$. As usual, we use the notation
$\ord X$ to denote the degree of the monomial $X$.

Let $1 \le i,j \le p, \; 1 \le l \le r_{ij}$ and let $n \ge N$. Then
we agree to write $f_{ij}(l)^{n}$ for the truncated power series
    \[ f_{ij}(l)^n = \sum_{\ord X = N}^n a_{ij}^l(X) \cdot X. \]
Moreover, let $a_{n+1} = I^{n+1} + (f^n)$ for all $n \ge N$, where
$(f^n) \subseteq \infdef$ is the ideal generated by $\{ f_{ij}(l)^n:
1 \le i,j \le p, \; 1 \le l \le r_{ij} \}$, and let $a_n = I^n$ for
$2 \le n \le N$. We write $H_n = H / I(H)^n$ as usual, then $H_n =
\infdef / a_n$ for all $n \ge 2$, in accordance with the notation in
the proof of theorem \ref{t:exh}.

Recall that $H_2 = \infdef_2$ and that $\xi_2 \in \defm(H_2)$
denotes the universal deformation with the property that the couple
$(H_2,\xi_2)$ represents $\defm$ restricted to $\acat p(2)$. We
have assumed that $a \subseteq I^N$, and this means that there exists a
lifting of $\xi_2$ to $H_N = \infdef / a_N = \infdef_N$. Let us
proceed to find such a lifting $M_N \in \defm(H_N)$ explicitly.

We choose to describe the deformation $M_N$ in terms of M-free
resolutions. Let us define $\ol B(N-1)$ to be the set of all
monomials in $\infdef$ of degree at most $N-1$. Then $\{ \ol X: X
\in \ol B(N-1) \}$ is a monomial basis of $H_N$, and any M-free
resolution of $M_N$ can be described by a family $\{ \alpha(X): X
\in \ol B(N-1) \}$ of $1$-cochains satisfying the following
conditions:
\begin{itemize}
\item $\alpha(e_i) = d_{\cpd i}$ for $1 \le i \le p$,
\item $\alpha(x_{ij}(l))$ is a $1$-cocycle representing
$x_{ij}(l)^*$ for $1 \le i,j \le p, \; 1 \le l \le d_{ij}$,
\item For all $Z \in \ol B(N-1)$ and for all $m \ge 0$, we have
\begin{equation*}
\sum_{\substack{X', X'' \in \ol B(N-1) \\ X' X'' = Z}}
\alpha(X'')_{m+1} \; \alpha(X')_m = 0.
\end{equation*}
\end{itemize}
We know that a family of $1$-cochains with the above properties
exists, since we can find a lifting $M_N$ of $\xi_2$ to $H_N$ and
this deformation must have some M-free resolution. So we choose
one such family $\{ \alpha(X): X \in \ol B(N-1) \}$ and fix this
choice. This means that we have fixed a deformation $M_N \in
\defm(H_N)$ with an M-free resolution given by the corresponding
differentials. So $(H_N,M_N)$ is a pro-representing hull for
$\defm$ restricted to $\acat p(N)$.

\begin{lem}\label{l:oex}
Let $\pi: R \to S$ be any small surjection in $\acat p$, let $\phi:
H \to S$ be any morphism, and denote by $M_{\phi} \in \defm(S)$ the
deformation induced by $\phi$. Then we can lift $\phi$ to a
morphism $\ol \phi: \infdef \to R$ making the diagram
\[
\xymatrix{
\infdef \ar[d] \ar[r]^{\ol \phi} & R \ar[d]^{\pi} \\
H \ar[r]_{\phi} & S }
\]
commutative, and the obstruction $o(\pi,M_{\phi})$ for lifting
$M_{\phi}$ to $R$ is given by
    \[ o(\pi,M_{\phi}) = \sum_{i,j,l} \; y_{ij}(l)^* \otimes
    \ol \phi(f_{ij}(l)). \]
\end{lem}
\begin{proof}
By construction and functoriality, the obstruction $o(\pi,
M_{\phi})$ is given as the restriction of the composition $\ol \phi
\circ o$ to the $k$-linear subspace $(\ext^2_A(M_j,M_i)^*)
\subseteq \obstrdef$. Since $\{ y_{ij}(l) \}$ is a $k$-linear
basis for this subspace, we get the desired expression for the
obstruction.
\end{proof}


Let us define $b_N \subseteq \infdef$ to be the ideal $b_N = I a_N +
a_N I = I^{N+1}$, and consider the natural map $r_N: R_N \to H_N$,
where $R_N = \infdef / b_N = \infdef_{N+1}$. By construction, $r_N$
is a small surjection in $\acat p$, and the natural surjection $\ol
\phi_N: \infdef \to R_N$ makes the diagram
\[
\xymatrix{
\obstrdef \ar[r]^o & \infdef \ar[d] \ar[rr]^{\ol \phi_N} & &
R_N \ar[d]^{r_N} \\
& H \ar[rr]_{\phi_N} & & H_N }
\]
commutative. Let $B'(N)$ be the set of all monomials in $\infdef$ of
degree $N$. Since $\ker(r_N) = I^N / I^{N+1}$, we see that $\{ \ol X:
X \in B'(N) \}$ is a monomial basis for $\ker(r_N)$. Moreover, let
$\ol{B'}(N) = B'(N) \cup \ol B(N-1)$. Then clearly $\{ \ol X: X \in
\ol{B'}(N) \}$ is a monomial basis for $R_N$.

Since $r_N$ is a small surjection, there is an obstruction $o(r_N,
M_N)$ for lifting $M_N$ to $R_N$, and we see from lemma \ref{l:oex}
that this obstruction can be expressed as
\begin{align*}
o(r_N,M_N) & = \sum_{i,j,l} \; y_{ij}(l)^* \otimes \ol
\phi_N(f_{ij}(l)) \\
& = \sum_{i,j,l} \; y_{ij}(l)^* \otimes \ol f_{ij}(l) \\
& = \sum_{i,j,l} \; \sum_{X \in B'(N)} y_{ij}(l)^* \otimes
(a_{ij}^l(X) \cdot \ol X),
\end{align*}
where $\ol f_{ij}(l)$ and $\ol X$ denote the images of $f_{ij}(l)$
and $X$ in $R_N$.

We say that the family $D(N) = \{ \alpha(X): X \in \ol B(N-1) \}$
of $1$-cochains is a \emph{defining system} for the matric Massey
products of order $N$,
    \[ \langle \ul x^* ; X \rangle \text{ for } X \in B'(N).
    \]
Let $X \in B'(N)$ be any monomial of type $(i,j)$. We define the
\emph{matric Massey product} $\langle \ul x^* ; X \rangle$ to be
the coefficient of $\ol X$ in the obstruction $o(r_N,M_N)$ above.
Then we immediately see that this matric Massey product has value
    \[ \langle \ul x^* ; X \rangle = \sum_{l=1}^{r_{ij}} \;
    a_{ij}^l(X) \cdot y_{ij}(l)^*. \]
In other words, the coefficient of $X$ in the power series
$f_{ij}(l)$ is given by the matric Massey product $\langle \ul
x^* ; X \rangle$ above as
    \[ a_{ij}^l(X) = y_{ij}(l)(\langle \ul x^* ; X \rangle) \]
for $1 \le l \le r_{ij}$.

We notice that the matric Massey products of order $N$ defined
above are \emph{immediately defined}. In other words, they can be
expressed in terms of the matric Massey products of section
\ref{s:mmp}. In fact, the defining system $D(N)$ induces a
defining system $\{ \alpha(X'): X' \dvd X, \; X' \neq X \}$ in
the sense of section \ref{s:mmp}, and the value of the
corresponding matric Massey product $\langle \ul x^* ; X \rangle$
is exactly the coefficient of $\ol X$ in the obstruction
$o(r_N,M_N)$.

On the other hand, we can calculate the obstruction $o(r_N,M_N)$
using the defining system $D(N)$, and therefore also the
coefficient of $\ol X$ in this obstruction for each $X \in B'(N)$.
A straight-forward calculation show that this coefficient is
given by the $2$-cocycle $y(X)$ defined by
\begin{equation*}
y(X)_m \; = \sum_{\substack{X', X'' \in \ol B(N-1) \\ X' X'' = X}}
\alpha(X'')_{m+1} \; \alpha(X')_m
\end{equation*}
for all $m \ge 0$. This means that the matric Massey product
$\langle \ul x^* ; X \rangle$ is represented by $y(X)$, so we can
easily calculate all matric Massey products of order $N$ using
the defining system $D(N)$. This determines the truncated power
series $f_{ij}(l)^N$, since we have
    \[ f_{ij}(l)^N = \sum_{X \in B'(N)} a_{ij}^l(X) \cdot X =
    \sum_{X \in B'(N)} y_{ij}(l)(\langle \ul x^* ; X \rangle)
    \cdot X \]
for $1 \le i,j \le p, \; 1 \le l \le r_{ij}$.

Let $h_N: H_{N+1} \to H_N$ be the natural map. Then $\ker(h_N) =
I^N / a_{N+1}$, so we can find a subset $B(N) \subseteq B'(N)$ of
monomials in $\infdef$ of degree $N$ such that $\{ \ol X: X \in
B(N) \}$ is a monomial basis for $\ker(h_N)$. Let $\ol B(N) = B(N)
\cup \ol B(N-1)$, then clearly $\{ \ol X: X \in \ol B(N) \}$ is a
monomial basis for $H_{N+1}$. So for each monomial $X \in \infdef$
with $\ord X \le N$, we have a unique relation in $H_{N+1}$ of the
form
    \[ \ol X = \sum_{X' \in \ol B(N)} \beta(X,X') \; \ol{X'},
    \]
with $\beta(X,X') \in k$ for all $X' \in \ol B(N)$. Since we have
$o(h_N,M_N)=0$, we deduce that
    \[ \sum_{\ord X = N} \langle x^* ; X \rangle \;
    \beta(X,X') = 0 \]
for all $X' \in B(N)$. Notice that $\beta(X,X')=0$ if the monomials
$X$ and $X'$ do not have the same type. Therefore, it makes sense
to consider the $1$-cocycle
    \[ \sum_{\ord X = N} \beta(X,X') \; y(X), \]
and by the relation above, this is a $1$-coboundary. It follows
that we can find a $1$-cochain $\alpha(X')$ such that
    \[ d \; \alpha(X') = - \sum_{\ord X = N} \beta(X,X') \;
    y(X), \]
and we fix such a choice. Consider the family $\{ \alpha(X): X \in
\ol B(N) \}$. This defines an M-free complex over $H_{N+1}$ if and
only if we have
    \[ \sum_{\ord X = N} \beta(X,Z) \sum_{\substack{X',X'' \in
    \ol B(N) \\ X' X'' = X}} \alpha(X'') \; \alpha(X') \; = \;
    0 \]
for all $Z \in \ol B(N)$. By the definition of $\alpha(X')$ when
$X' \in B(N)$, this condition holds, and we denote by $M_{N+1} \in
\defm(H_{N+1})$ the deformation with the complex defined by $\{
\alpha(X): X \in \ol B(N) \}$ as M-free resolution. It is clear
from the construction that $M_{N+1}$ is a lifting of $M_N$, so
$(H_{N+1},M_{N+1})$ is a pro-representing hull for $\defm$
restricted to $\acat p(N+1)$.


Let $b_{N+1} \subseteq \infdef$ be the ideal $b_{N+1} = I a_{N+1}
+ a_{N+1} I = I^{N+2} + I (f^N) + (f^N) I$, and consider the
natural map $r_{N+1}: R_{N+1} \to H_{N+1}$, where $R_{N+1} =
\infdef / b_{N+1}$. By construction, $r_{N+1}$ is a small
surjection in $\acat p$, and it is clear that the natural
morphism $\ol \phi_{N+1}: \infdef \to R_{N+1}$ makes the diagram
\[
\xymatrix{
\obstrdef \ar[r]^o & \infdef \ar[d] \ar[rr]^{\ol \phi_{N+1}} & &
R_{N+1} \ar[d]^{r_{N+1}} \\
& H \ar[rr]^{\phi_{N+1}} \ar[rrd]_{\phi_N} & & H_{N+1}
\ar[d]^{h_N} \\
& & & H_N }
\]
commutative. We see that $\ker(r_{N+1}) = a_{N+1}/b_{N+1}$, which
we can re-write in the following way:
\begin{align*}
\ker(r_{N+1}) & = ( I^{N+1} + (f^N) ) / ( I^{N+2} + I(f^N) + (f^N)
I ) \\ & = (f^N) / ( I (f^N) + (f^N) I ) \; \oplus \; I^{N+1} /
( I^{N+2} + I (f^N) + (f^N) I )
\end{align*}
Let us write $c(N+1) = I^{N+1} / ( I^{N+2} + I (f^N) + (f^N) I )$.
Then $c(N+1) \subseteq \ker(r_{N+1})$ is an ideal, and we can
clearly find a set $B'(N+1)$ of monomials in $\infdef$ of degree
$N+1$ such that $\{ \ol X: X \in B'(N+1) \}$ is a monomial basis
for $c_{N+1}$. Let us choose $B'(N+1)$ such that for every $X \in
B'(N+1)$, there is a monomial $X' \in B(N)$ such that $X' \dvd X$,
this is clearly possible. We let $\ol{B'}(N+1) = B'(N+1) \cup \ol
B(N)$, then
    \[ \{ \ol X: X \in \ol{B'}(N+1) \} \cup \{ f_{ij}(l)^N :
    1 \le i,j \le p, \; 1 \le l \le r_{ij} \} \]
is a basis for $R_{N+1}$. So for each monomial $X \in \infdef$
with $\ord X \le N+1$, we have a unique relation in $R_{N+1}$ of
the form
    \[ \ol X = \sum_{X' \in \ol{B'}(N+1)} \beta'(X,X')
    \; \ol{X'} \; + \; \sum_{i,j,l} \; \beta'(X,i,j,l) \;
    \ol f_{ij}(l)^N, \]
with $\beta'(X,X'), \beta'(X,i,j,l) \in k$ for all $X' \in \ol{B'}
(N+1), \; 1 \le i,j \le p, \; 1 \le l \le r_{ij}$.

Since $r_{N+1}$ is a small surjection, there is an obstruction
$o(r_{N+1},M_{N+1})$ for lifting $M_{N+1}$ to $R_{N+1}$, and we
see from lemma \ref{l:oex} that this obstruction can be expressed
as
\begin{align*}
o(r_{N+1},M_{N+1}) & = \sum_{i,j,l} \; y_{ij}(l)^* \otimes \ol
\phi_{N+1}(f_{ij}(l)) \\
& = \sum_{i,j,l} \; y_{ij}(l)^* \otimes \ol f_{ij}(l) \\
& = \sum_{i,j,l} \; y_{ij}(l)^* \otimes ( \ol f_{ij}(l)^N +
\sum_{X \in B'(N+1)} a_{ij}^l(X) \cdot \ol X ),
\end{align*}
where $\ol f_{ij}(l), \; \ol f_{ij}(l)^N$ and $\ol X$ denote the
images of $f_{ij}(l), \; f_{ij}(l)^N$ and $X$ in $R_{N+1}$.

We say that the family $D(N+1) = \{ \alpha(X): X \in \ol B(N)
\}$ is a \emph{defining system} for the matric Massey products
of order $N+1$,
    \[ \langle \ul x^* ; X \rangle \text{ for } X \in
    B'(N+1) \]
Let $X \in B'(N+1)$ be any monomial of type $(i,j)$. We define
the \emph{matric Massey product} $\langle \ul x^* ; X \rangle$
to be the coefficient of $\ol X$ in the obstruction $o(r_{N+1},
M_{N+1})$ above. Then we immediately see that this matric Massey
product has value
    \[ \langle \ul x^* ; X \rangle = \sum_{l=1}^{r_{ij}} \;
    a_{ij}^l(X) \cdot y_{ij}(l)^*. \]
In other words, the coefficient of $X$ in the power series
$f_{ij}(l)$ is given by the matric Massey product $\langle \ul
x^* ; X \rangle$ above as
    \[ a_{ij}^l(X) = y_{ij}(l)(\langle \ul x^* ; X \rangle)
    \]
for $1 \le l \le r_{ij}$.

On the other hand, we can calculate the obstruction $o(r_{N+1},
M_{N+1})$ using the defining system $D(N+1)$, and therefore also
the coefficient of $\ol X$ in this obstruction for each $X \in
B'(N+1)$. A straight-forward calculation show that this
coefficient is given by the $2$-cocycle $y(X)$ defined by
\begin{equation*}
y(X)_m \; = \sum_{\ord Z \le N+1} \beta'(Z,X)
\sum_{\substack{X', X'' \in \ol B(N) \\ X' X'' = Z}}
\alpha(X'')_{m+1} \; \alpha(X')_m
\end{equation*}
for all $m \ge 0$. This means that the matric Massey product
$\langle \ul x^* ; X \rangle$ is represented by $y(X)$, so we can
easily calculate all matric Massey products of order $N+1$ using
the defining system $D(N+1)$.

By the construction in the proof of theorem \ref{t:exh}, we have
that $H_{N+2}$ is the quotient of $R_{N+1}$ by the ideal generated
by the obstruction $o(r_{N+1},M_{N+1})$. On the other hand, we
know that $H_{N+2} = \infdef / (I^{N+2} + (f^{N+1})$. This implies
that for all monomials $X \not \in B'(N+1)$ of degree $N+1$, the
coefficient $a_{ij}^l(X) = 0$ for $1 \le i,j \le p, \; 1 \le l \le
r_{ij}$. In other words, the truncated power series
$f_{ij}(l)^{N+1}$ is determined by the matric Massey products of
order $N+1$ above, since we have
\begin{align*}
f_{ij}(l)^{N+1} & = f_{ij}(l)^N + \sum_{X \in B'(N+1)}
a_{ij}^l(X) \cdot X \\
& = f_{ij}(l)^N + \sum_{X \in B'(N+1)} y_{ij}(l)(\langle
\ul x^* ; X \rangle) \cdot X
\end{align*}
for $1 \le i,j \le p, \; 1 \le l \le r_{ij}$.

Let $h_{N+1}: H_{N+2} \to H_{N+1}$ be the natural map, and
consider its kernel. By definition, we have
    \[ \ker(h_{N+1}) = a_{N+1}/a_{N+2} = ((f^N) + I^{N+1})/
    ((f^{N+1}) + I^{N+2}), \]
so we can clearly find a subset $B(N+1) \subseteq B'(N+1)$ of
monomials of degree $N+1$ such that $\{ \ol X: X \in \ol B(N+1) \}
\cup \{ \ol f_{ij}(l)^N: 1 \le i,j \le p, \; 1 \le l \le r_{ij}
\}$ is a basis for $\ker(h_{N+1})$. Let $\ol B(N+1) = B(N+1) \cup
\ol B(N)$, then clearly
    \[ \{ \ol X: X \in \ol B(N+1) \} \cup \{ \ol f_{ij}(l)^N:
    1 \le i,j \le p, \; 1 \le l \le r_{ij} \} \]
is a monomial basis for $H_{N+2}$. So for each monomial $X \in
\infdef$ with $\ord X \le N+1$, we have a unique relation in
$H_{N+2}$ of the form
    \[ \ol X = \sum_{X' \in \ol B(N+1)} \beta(X,X') \;
    \ol{X'} \; + \; \sum_{i,j,l} \; \beta(X,i,j,l) \;
    f_{ij}(l)^N, \]
with $\beta(X,X'), \beta(X,i,j,l) \in k$ for all $X' \in \ol
B(N+1), \; 1 \le i,j \le p, \; 1 \le l \le r_{ij}$. Since we
have $o(h_{N+1},M_{N+1})=0$, we deduce that
    \[ \sum_{\ord X \le N+1} \langle x^* ; X \rangle \;
    \beta(X,X') = 0 \]
for all $X' \in B(N+1)$. Notice that $\beta(X,X')=0$ if the
monomials $X$ and $X'$ do not have the same type. Therefore, it
makes sense to consider the $1$-cocycle
    \[ \sum_{\ord X \le N+1} \beta(X,X') \; y(X), \]
and by the relation above, this is a $1$-coboundary. It follows
that we can find a $1$-cochain $\alpha(X')$ such that
    \[ d \; \alpha(X') = - \sum_{\ord X \le N+1} \beta(X,X')
    \; y(X), \]
and we fix such a choice. Consider the family $\{ \alpha(X): X \in
\ol B(N+1) \}$. This defines an M-free complex over $H_{N+2}$ if
and only if we have
    \[ \sum_{\ord X \le N+1} \beta(X,Z) \sum_{\substack{X',X''
    \in \ol B(N+1) \\ X' X'' = X}} \alpha(X'') \; \alpha(X') \;
    = \; 0 \]
for all $Z \in \ol B(N+1)$. By the definition of $\alpha(X')$ when
$X' \in B(N+1)$, this condition holds, and we denote by $M_{N+2}
\in \defm(H_{N+2})$ the deformation with the complex defined by
$\{ \alpha(X): X \in \ol B(N+1) \}$ as M-free resolution. It is
clear from the construction that $M_{N+2}$ is a lifting of
$M_{N+1}$, so $(H_{N+2},M_{N+2})$ is a pro-representing hull for
$\defm$ restricted to $\acat p(N+2)$.

It is clear that we can continue in this way. For every $k \ge 1$,
we can calculate the coefficients in the truncated power series
$f_{ij}(l)^{N+k}$, and therefore find $H_{N+k+1}$. At the same
time, we find the defining systems $\{ \alpha(X): X \in \ol B(N+k)
\}$ necessary to calculate the matric Massey products of order
$N+k+1$, and these defining systems completely determine the
deformation $M_{N+k+1}$. We have described how to do this in the
case $k=1$, and the general case is similar.

We conclude that the method that we have described above can be
used to calculate the pro-representing hull $(H_n,M_n)$ for the
deformation functor $\defm$ restricted to $\acat p(n)$ for any $n
\ge N$. We can therefore, in principle, find the hull
    \[ H = \lim_{\leftarrow} H_n \]
of $\defm$, and also the corresponding versal family defined over
$H$,
    \[ \xi = M = \lim_{\leftarrow} M_n. \]
It follows that the pro-representing hull $(H,\xi)$ of the
deformation functor $\defm$ can be calculated using matric Massey
products.

\section{An example}

Let $k$ be an algebraically closed field of characteristic $0$,
and let $A = A_2(k)$ be the second Weyl algebra over $k$. We shall
think of $A$ as the ring of differential operators in the plane
defined over $k$ with coordinates $x$ and $y$. Thus, we can write
$A = k[x,y] \langle \dx, \dy \rangle$, where $\dx = \partial /
\partial x$ and $\dy = \partial / \partial y$. In other words, $A$
is the $k$-algebra generated by $x,y,\dx,\dy$ with relations
$[\dx, x] = [\dy, y] = 1$.

Let us consider the family of left $A$-modules $\mfam = \{ M_1,
M_2, M_3, M_4 \}$, where $M_i = A/I_i$ for $1 \le i \le 4$ and
$I_i \subseteq A$ are left ideals given by
\begin{align*}
    I_1 &= A ( \dx, \dy )   & I_2 &= A ( \dx, y ) \\
    I_3 &= A ( x, \dy)    & I_4 &= A ( x, y )
\end{align*}
We immediately notice that the left $A$-modules in the family
$\mfam$ have the following free resolutions:
\begin{align*}
0 \gets &M_1 \gets A \xleftarrow{\begin{pmatrix} \dx \\ \dy
\end{pmatrix}} A^2 \xleftarrow{\begin{pmatrix} \dy & -\dx
\end{pmatrix}} A \gets 0 \\
0 \gets &M_2 \gets A \xleftarrow{\begin{pmatrix} \dx \\ y
\end{pmatrix}} A^2 \xleftarrow{\begin{pmatrix} y & -\dx
\end{pmatrix}} A \gets 0 \\
0 \gets &M_3 \gets A \xleftarrow{\begin{pmatrix} x \\ \dy
\end{pmatrix}} A^2 \xleftarrow{\begin{pmatrix} \dy & -x
\end{pmatrix}} A \gets 0 \\
0 \gets &M_4 \gets A \xleftarrow{\begin{pmatrix} x \\ y
\end{pmatrix}} A^2 \xleftarrow{\begin{pmatrix} y & -x
\end{pmatrix}} A \gets 0 \\
\end{align*}
We consider the elements of the free $A$-modules $A^n$ as row
vectors, and the maps in the free resolutions above as right
multiplication of these row vectors by the given matrices. Notice
that for $1 \le i \le 4$, the free $A$-module $L_{m,i}$ in the
free resolution of $M_i$ does not depend upon $i$. We shall
therefore write $L_m = L_{m,i}$ for all $m \ge 0, \; 1 \le i \le
4$.

It is known that $\mfam$ is a family of simple holonomic left
$A$-modules, so this family satisfy the finiteness condition
(\ref{e:fc}). Therefore, there exists a pro-representing hull
$(H,\xi)$ for the deformation functor $\defm: \acat 4 \to \sets$
by theorem \ref{t:exh}. We shall use the methods from section
\ref{s:ch} to construct this hull explicitly.

Let us start by calculating $\ext^n_A(M_i,M_j)$ for $n=1,2, \; 1
\le i,j \le 4$. We need both the dimensions and $k$-linear bases
for these vector spaces, where each basis vector is represented by
a cocycle in the corresponding Yoneda complex. The calculations
are straight-forward, so we only state the results here:
\begin{align*}
\dim_k \ext^1_A(M_i,M_j) &=
\begin{cases}
    1& \text{if $i=1$ or $i=4$ and $j=2$ or $j=3$, or} \\
    & \text{if $i=2$ or $i=3$ and $j=1$ or $j=4$,} \\
    0& \text{otherwise}
\end{cases} \\
\intertext{}
\dim_k \ext^2_A(M_i,M_j) &=
\begin{cases}
    1& \text{if $(i,j) = (1,4), (2,3), (3,2), (4,1)$,} \\
    0& \text{otherwise}
\end{cases}
\end{align*}
We denote the basis vectors of $\ext^1_A(M_j,M_i)$ by $x_{ij}^*$
since there is at most one for each pair of indices $(i,j)$. From
the dimensions listed above, we see that we have the following
basis vectors:
    \[ x_{12}^*, x_{13}^*, x_{21}^*, x_{24}^*, x_{31}^*, x_{34}^*,
    x_{42}^*, x_{43}^* \]
We choose a Yoneda representative for each vector $x_{ij}^*$ in
this list, and we denote this representative by $\alpha(x_{ij})$.
From the free resolutions above, we see that we can write each of
these representatives in the form
    \[ \alpha(X) = \{ \alpha(X)_0, \alpha(X)_1 \}, \]
where $\alpha(X)_0: L_1 \to L_0$ is right multiplication by a
matrix $\left(\begin{smallmatrix} a \\ b \end{smallmatrix}\right)$
with entries $a,b \in A$, and $\alpha(X)_1: L_2 \to L_1$ is right
multiplication by a matrix $\left(\begin{smallmatrix} c & d
\end{smallmatrix}\right)$ with entries $c,d \in A$ for each
monomial $X = x_{ij}$. We find the following representatives:
\begin{align*}
\alpha(x_{12}) = \alpha(x_{21}) = \alpha(x_{34}) = \alpha(x_{43})
&= \{ \left( \begin{smallmatrix} 0 \\ 1 \end{smallmatrix} \right),
\; \left( \begin{smallmatrix} 1 & 0 \end{smallmatrix} \right) \}
\\
\alpha(x_{13}) = \alpha(x_{31}) = \alpha(x_{24}) = \alpha(x_{42})
&= \{ \left( \begin{smallmatrix} 1 \\ 0 \end{smallmatrix} \right),
\; \left( \begin{smallmatrix} 0 & -1 \end{smallmatrix} \right) \}
\end{align*}
Similarly, we denote the basis vectors of $\ext^2_A(M_j,M_i)$ by
$y_{ij}^*$ since there is at most one for each pair of indices
$(i,j)$. From the dimensions listed above, we see that we have the
following basis vectors:
    \[ y_{14}^*, y_{23}^*, y_{32}^*, y_{41}^* \]
We choose a Yoneda representative for each vector $y_{ij}^*$ in
this list, and we denote this representative $\alpha(y_{ij})$.
From the free resolutions above, we see that we can write each of
these representatives in the form
    \[ \alpha(Y) = \{ \alpha(Y)_0 \}, \]
where $\alpha(Y)_0: L_2 \to L_0$ is given by right multiplication
of an element $a \in A$ for each monomial $Y = y_{ij}$. We find
the following representatives:
\begin{align*}
\alpha(y_{14}) = \alpha(y_{23}) = \alpha(x_{32}) = \alpha(x_{41})
&= \{ \left( \begin{smallmatrix} 1 \end{smallmatrix} \right) \}
\end{align*}
This completes the calculations of $\ext^n_A(M_i,M_j)$ for $n=1,2$
and $1 \le i,j \le 4$. We know that these calculations determine
the hull at the tangent level, $(H_2,\xi_2)$.

The next step is to find the the hull $H$ and the versal family
$\xi$, and we shall employ the notations and methods of section
\ref{s:ch} to accomplish this. Let $N = 2$, we know that this
choice is always possible. As usual, we let $\infdef$ be the
formal matrix algebra generated by the monomials $x_{ij}$ in the
above list, and let $I = I(\infdef)$ be its radical. Furthermore,
denote by $f_{ij} = o(y_{ij}) \in I^2_{ij}$ for $(i,j) = (1,4),
(2,3), (3,2), (4,1)$, and by $f_{ij}^n$ the corresponding
truncated power series for each $n \ge N$.

First, we have to find a defining system $\{ \alpha(X): \ord X < 2
\}$ for the matric Massey products $\langle \ul x^* ; X \rangle$
when $X$ is any monomial of degree $2$ in $\infdef$. This is
easily done: The $1$-cocycle $\alpha(e_i)$ is the free resolution
of $M_i$ for $1 \le i \le 4$, and the $1$-cocycle $\alpha(X)$ was
chosen above for each monomial $X = x_{ij}$ of degree $1$.

Let us calculate the matric Massey products of order $2$: Using
the defining system given above, we find that the cocycles $y(X)$
representing the matric Massey products $\langle \ul x^* ; X
\rangle$ are given by
\[
y(X)_0 =
\begin{cases}
    -1& \text{if $X = x_{12} x_{24}, x_{21} x_{13}, x_{34} x_{42},
    x_{43} x_{31}$,} \\
    1&  \text{if $X = x_{13} x_{34}, x_{24} x_{43}, x_{31} x_{12},
    x_{42} x_{21}$,} \\
    0& \text{otherwise}
\end{cases}
\]
for all monomials $X$ of degree $2$ in $\infdef$. This means that
the corresponding matric Massey products are given by
\begin{align*}
    \langle x_{12}, x_{24} \rangle &= - y_{14}   &
    \langle x_{13}, x_{34} \rangle &= y_{14} \\
    \langle x_{21}, x_{13} \rangle &= - y_{23}   &
    \langle x_{24}, x_{43} \rangle &= y_{23} \\
    \langle x_{31}, x_{12} \rangle &= y_{32}   &
    \langle x_{34}, x_{42} \rangle &= - y_{32} \\
    \langle x_{42}, x_{21} \rangle &= y_{41}   &
    \langle x_{43}, x_{31} \rangle &= - y_{41},
\end{align*}
and all other matric Massey products of order $2$ are zero. This
translates to the following truncated power series $f_{ij}^2$:
\begin{align*}
    f_{14}^2 &= x_{13} x_{34} - x_{12} x_{24} \\
    f_{23}^2 &= x_{24} x_{43} - x_{21} x_{13} \\
    f_{32}^2 &= x_{31} x_{12} - x_{34} x_{42} \\
    f_{41}^2 &= x_{42} x_{21} - x_{43} x_{31}
\end{align*}
By the general theory, we therefore have $H_3 = \infdef /
(f_{14}^2, f_{23}^2, f_{32}^2, f_{41}^2) + I^3$. We know that we
can find a lifting $\xi_3$ of $\xi_2$ to $H_3$, and that
$(H_3,\xi_3)$ is a pro-representing hull of $\defm$ restricted to
$\acat 4(3)$.

In order to find $\xi_3$, we let $B(2) = \{ X: \ord X = 2 \}
\setminus \{ x_{13} x_{34}, x_{24} x_{43}, x_{31} x_{12}, x_{42}
x_{21} \}$. We also let $\ol B(2) = B(2) \cup \ol B(1)$, where
$\ol B(1) = \{ X: \ord X \le 1 \}$. Then $\{ \ol X: X \in \ol B(2)
\}$ is a monomial basis for $H_3$. We observe that if we choose
$\alpha(X) = 0$ for all $X \in B(2)$, the family $\{ \alpha(X): X
\in \ol B(2) \}$ defines an M-free complex over $H_3$. In other
words, this family completely defines the deformation $\xi_3 \in
\defm(H_3)$ lifting $\xi_2$.

Clearly, we could continue in this way. But after the last
computations, it is tempting to think that $f_{ij} = f_{ij}^2$ for
$(i,j) = (1,4), (2,3), (3,2), (4,1)$. Let us check if this is the
case: We put $T = \infdef / (f_{14}^2, f_{23}^2, f_{32}^2,
f_{41}^2)$, and choose a monomial basis $B$ of $T$ containing $\ol
B(2)$. Furthermore, we let $\alpha(X)$ be as before when $X \in
\ol B(2)$ and let $\alpha(X) = 0$ for all monomials $X \in B$ of
degree at least $3$. This choice corresponds to maps $d^T_0,
d^T_1$ of M-free modules over $T$, and a computation shows that
    \[ d^T_0 \circ d^T_1 = (1 \otimes (f_{14}^2 + f_{23}^2 +
    f_{32}^2 + f_{41}^2)) = 0. \]
So the family $\{ \alpha(X): X \in B \}$ defines an M-free complex
over $T$, and therefore a deformation $\xi \in \defm(T)$ lifting
$\xi_3$. This proves that $H=T$, or in other words, that
    \[ H = \infdef / (x_{13} x_{34} - x_{12} x_{24}, x_{24}
    x_{43} - x_{21} x_{13}, x_{31} x_{12} - x_{34} x_{42}, x_{42}
    x_{21} - x_{43} x_{31}) \]
is a pro-representing hull of $\defm$. In particular, $f_{ij} =
f_{ij}^2$ for all $i,j$. Moreover, the family $\{ \alpha(X): X \in
B \}$ defines the versal family $\xi \in
\defm(H)$.

\appendix

\section{Yoneda and Hochschild representations}
\label{s:app} \label{a:hh}

Let $k$ be an algebraically closed (commutative) field, let $A$ be
an associative $k$-algebra, and let $M,N$ be left $A$-modules. In
this appendix, we recall several different descriptions of the
$k$-vector space $\ext_A^n(M,N)$ for $n \ge 0$. In particular, we
show how to realize this cohomology group using the Yoneda and
Hochschild complexes.

\subsection{The Yoneda representation}

Fix free resolutions $(L_{\cpd},d_{\cpd})$ of $M$ and $(L'_\cpd,d'_\cpd)$ of $N$.
We shall write $d_i: L_{i+1} \to L_i$ and $d'_i: L'_{i+1} \to L'_i$ for the
differentials, and denote the augmentation morphisms by $\rho: L_0 \to M$ and
$\rho': L'_0 \to N$.

For all integers $n\ge 0$, the cohomology group $\ext^n_A(M,N)$ is
defined to be the $n$'th cohomology group of the complex
$Hom_A(L_\cpd,N)$,
\[ \ext_A^n(M,N) = H^n(Hom_A(L_\cpd,N)). \]
Notice that in general, this Abelian group does not have a left
$A$-module structure, but only a left $\cen(A)$-module structure,
where $\cen(A)$ is the centre of $A$. In particular, if $A$ is
commutative, then $\ext_A^n(M,N)$ has the structure of an
$A$-module, and if $A$ is a $k$-algebra, then $\ext^n_A(M,N)$ has
the structure of a $k$-vector space.

We denote by $\hmm^\cpd(L_\cpd,L'_\cpd)$ the \emph{Yoneda complex}
given by the given free resolutions. This complex is defined in
the following way: For each integer $n\ge 0$, let
$\hmm^n(L_\cpd,L'_\cpd)$ be the left $A$-module
\[ \hmm^n(L_\cpd,L'_\cpd) = \amalg_i \hmm_A(L_{i+n},L'_i). \]
Moreover, let the differential $d^n: \hmm^n(L_\cpd,L'_\cpd) \to
\hmm^{n+1}(L_\cpd, L'_\cpd)$ for $n \ge 0$ be the $A$-linear map
given by the formula
\[ d^n(\phi)_i = \phi_i d_{n+i} + (-1)^{n+1} d'_i \phi_{i+1} \]
for all $i\ge 0$, where we write $\phi = (\phi_i)$ with $\phi_i
\in \hmm_A(L_{i+n}, L'_i)$ for all $i\ge 0$. It is easy to check
that this map is a well-defined differential, so the Yoneda
complex is a complex of Abelian groups. We shall write
$H^n(\hmm(L_\cpd,L'_\cpd))$ for the cohomology groups of the
Yoneda complex. Since the differential $d=d^n$ is left
$\cen(A)$-linear, these cohomology groups have a natural structure
as left $\cen(A)$-modules.

\begin{lem}
\label{l:yonext} For all integers $n\ge 0$, there is a canonical
isomorphism of left $\cen(A)$-modules
\[ H^n(\hmm(L_\cpd,L'_\cpd)) \cong \ext_A^n(M,N). \]
\end{lem}
\begin{proof} \noindent
There is a natural map $f_n: \hmm^n(L_\cpd,L'_\cpd) \to
\hmm_A(L_n,N)$, given by $f(\phi)=\rho' \phi_0$, where $\phi =
(\phi_i) \in \hmm^n(L_\cpd,L'_\cpd)$. It is easy to see that these
maps are compatible with the differentials, and a small
calculation show that $f_n$ induces an isomorphism on cohomology
$H^n(\hmm(L_\cpd, L'_\cpd)) \to \ext_A^n(M,N)$ for all integers
$n\ge 0$.
\end{proof}

\subsection{Definition of Hochschild cohomology}

Let $Q$ be an $A$-$A$ bimodule. We define the \emph{Hochschild
complex} of $A$ with values in $Q$ in the following way: Let
$\hc^n(A,Q) = \hmm_k(\otimes^n_k A,Q)$ for all $n\ge 0$. So any
$\psi \in \hc^n(A,Q)$ corresponds to a $k$-multilinear map from
$n$ copies of $A$ into $Q$, and we shall therefore write
$\psi(a_1, \dots, a_n)$ in place of $\psi(a_1 \otimes \dots
\otimes a_n)$ for $\psi\in \hc^n(A,Q), \; a_1, \dots, a_n \in A$.
Moreover, let $d^n: \hc^n(A,Q) \to \hc^{n+1}(A,Q)$ for $n \ge 0$
be the $k$-linear map given by the formula
\begin{multline}
\label{f:hhdf} d^n (\psi)(a_0, \dots , a_n) = a_0 \, \psi(a_1,
\dots ,a_n) +
\sum_{i=1}^n (-1)^i \psi(a_0, \dots ,a_{i-1}a_i, \dots ,a_n) \\
+ (-1)^{n+1} \psi(a_0, \dots ,a_{n-1}) \, a_n
\end{multline}
for all $\psi\in \hc^n(A,Q), \; a_0, \dots, a_n \in A$.

\begin{lem}
$\hc^*(A,Q)$ is a complex of $k$-vector spaces.
\end{lem}
\begin{proof} \noindent
Let $\psi \in \hc^n(A,Q)$. Then $\psi' = d^n(\psi)$ is a sum of
$n+1$ summands, and we denote these by $\psi'_0, \dots , \psi'_n$,
in the order they appear in formula \ref{f:hhdf}. We let $\psi'' =
d^{n+1} \psi' = d^{n+1} d^n \psi$. Each $d^{n+1} \psi'_i$ for $0
\le i \le n$ is a sum of $n+2$ summands, and we denote these by
$\psi''_{ij}$ for $0 \le j \le n+1$ in the order they appear in
formula \ref{f:hhdf}. A straight-forward calculation shows that we
have $\psi''_{i,j} + \psi''_{j,i+1} = 0$ for all indices $i,j$
with $0 \le j \le n+2, \; j \le i \le n+1$. Since $\psi'' = \sum
\psi''_{ij}$, it follows that $\psi''=0$ in $\hc^{n+2}(A,Q)$.
Consequently, $\hc^*(A,Q)$ is a complex of $k$-vector spaces.
\end{proof}

We define the \emph{Hochschild cohomology} of $A$ with values in
$Q$ to be the cohomology of the Hochschild complex
$\hc^\cpd(A,Q)$, so we have
\[ \hh^n(A,Q) = H^n(\hc^*(A,Q)) = \ker(d^n) / \im(d^{n-1}) \]
for all $n\ge 0$. In particular, the cohomology groups
$\hh^n(A,Q)$ have a natural structure as $k$-vector spaces.

Let $\psi \in \hc^1(A,Q)$, then $\psi$ is a $1$-cocycle if and
only if $\psi(ab) = a \psi(b) + \psi(a) b$ for all $a,b \in A$. So
we have $\ker(d^1)=\der_k(A,Q)$. We say that a derivation $\psi\in
\der_k(A,Q)$ is trivial if there is an element $q\in Q$ such that
$\psi$ is of the form $\psi(a)=aq-qa$ for all $a \in A$. Clearly,
the set of trivial derivations is the image $\im(d^0)$. So
$\hh^1(A,Q) \cong \der_k(A,Q)/T$ where $T$ is the trivial
derivations of $A$ into $Q$.

\subsection{The Hochschild representation}

We remark that $Q=\hmm_k(M,N)$ is an $A$-$A$ bimodule in a natural
way: For any $a \in A$, let $L_a: M \to M$ denote left
multiplication on $M$ by $a$, and $L'_a: N \to N$ left
multiplication on $N$ by $a$. The bimodule structure is given by
$a \phi = L'_a \phi, \; \phi a = \phi L_a$ for $a \in A, \; \phi
\in \hmm_k(M,N)$. We shall consider the Hochschild cohomology of
$A$ with values in $Q=\hmm_k(M,N)$.

By definition, we have that $\hh^0(A,Q) = \hmm_A(M,N)$ when
$Q=\hmm_k(M,N)$. So we have a natural isomorphism of $k$-vector
spaces $\ext^0_A(M,N) \cong \hh^0(A,Q)$. Notice that since $k
\subseteq \cen(A)$, $\ext_A^n(M,N)$ has a natural $k$-vector space
structure for all $n\ge 0$. It is possible to extend the above
isomorphism to the higher cohomology groups:

\begin{prop}
\label{p:hhext}
For all integers $n \ge 0$, there is an isomorphism of $k$-vector spaces
\[ \sigma_n : \ext^n_A(M,N) \to \hh^n(A,\hmm_k(M,N)). \]
\end{prop}
\begin{proof} \noindent
From Weibel \cite{wei94}, lemma 9.1.9, there is an isomorphism of
$k$-vector spaces between $\hh^n(A,\hmm_k(M,N))$ and
$\ext^n_{A/k}(M,N)$ for $n\ge 0$. But since $k$ is a commutative
field, there is a canonical isomorphism between
$\ext^n_{A/k}(M,N)$ and $\ext^n_A(M,N)$, see theorem 8.7.10 in
Weibel \cite{wei94}.
\end{proof}

We shall give an explicit identification of $k$-vector spaces
between $\ext_A^1(M,N)$ and $\hh^1(A,\hmm_k(M,N))$: Let
$(L_\cpd,d_\cpd)$ be a free resolution of $M$, with augmentation
morphism $\rho: L_0 \to M$, and let $\tau: M \to L_0$ be a
$k$-linear section of $\rho$. For any $1$-cocycle $\phi \in
\hmm_A(L_1,N)$, let $\psi = \psi(\phi) \in \der_k(A,\hmm_k(M,N))$
be the following derivation: For any $a \in A, \; m \in M$, let
$x=x(a,m) \in L_1$ be such that $d_0(x) = a \tau(m) - \tau(am)$.
Notice that such an $x$ exists, and is uniquely defined modulo the
image $\im d_1$. We define $\psi$ by the equation
$\psi(a)(m)=\phi(x)$ with $x=x(a,m)$. Since $\phi$ is a cocycle,
$\psi$ is a well-defined homomorphism in $\hmm_k(A,\hmm_k(M,N))$,
and a straight-forward calculation shows that $\psi$ is a
derivation.

\begin{lem}
Assume that $\ext_A^1(M,N)$ is a finite dimensional $k$-vector
space. Then the assignment $\phi \mapsto \psi(\phi)$ defined in
the above paragraph induces an isomorphism $\sigma_1:
\ext_A^1(M,N) \to \hh^1(A,\hmm_k(M,N))$.
\end{lem}
\begin{proof} \noindent
Assume that $\phi$ is a co-boundary, so $\phi=d^0(\phi')$, where
$\phi' \in \hmm_A(L_0,N)$. Then $\psi = d^0(\phi')$, where
$\psi'=\phi' \tau \in \hmm_k(M,N)$, so $\phi$ is a trivial
derivation. Consequently, the assignment induces a well-defined
map of $k$-linear spaces. This map is furthermore injective:
Assume that $\psi$ is a trivial derivation, so $\psi =
d^0(\psi')$, where $\psi' \in \hmm_k(M,N)$. Then, we can construct
an $A$-linear map $\phi' \in \hmm_A(L_0,N)$ in the following way:
Choose a basis for $L_0$, and for each basis vector $y\in L_0$,
choose $y' \in L_1$ such that $d_0(y')=y-\psi' \rho(y)$. Then we
define $\phi'(y) = \psi' \rho (y) + \phi(y')$ for each basis
vector $y\in L_0$. We obtain a morphism $\phi' \in \hmm_A(L_0,N)$
by $A$-linear extension, and $d^0(\phi')=\phi$, so $\phi$ is a
co-boundary. To show that $\sigma_1$ is an isomorphism as well, it
is enough to notice that $\dim_k \ext_A^1(M,N) = \dim_k
\hh^1(A,\hmm_k(M,N))$ by proposition \ref{p:hhext}, since
$\ext^1_A(M,N)$ has finite $k$-dimension.
\end{proof}

The identification $\sigma_n: \ext^n_A(M,N) \to
\hh^n(A,\hmm_k(M,N))$ for $n \ge 2$ can be constructed in a
similar way.

\bibliographystyle{amsplain}
\bibliography{commalg,noncommalg,homalg,deform,diff}

\end{document}